\documentclass[11pt]{amsart}

\usepackage[latin1]{inputenc}%para poner los tildes directamente
\usepackage{amssymb}%para las letras dobles
\usepackage{amsmath}
\usepackage{esint}
\usepackage{hyperref}
\usepackage{tikz, pgf, xcolor, graphicx}
\usepackage{fancyvrb}
\usepackage{color}
\usepackage{url}

\hypersetup{
    colorlinks=true,
    linkcolor=black,
    citecolor=black,
    filecolor=black,
    urlcolor=black,
}

\newcommand{\rn}{{\mathbb{R}^n}}
\newcommand{\phii}{\varphi}

\newcommand{\w}{\omega}
\newcommand{\W}{\Omega}
\newcommand{\x}{\hat {{x}}}
\newcommand{\y}{\hat {{y}}}
\newcommand{\z}{\hat {{z}}}
\newcommand{\supp}{\text{supp}}
\newcommand{\eps}{\varepsilon}
\usepackage{color}

\allowdisplaybreaks

\def\R{{\mathbb {R}}}
\def\N{{\mathbb {N}}}

\def\V{{\mathbb {V}}}

\newtheorem{theorem}{Theorem}[section]

\newtheorem{proposition}[theorem]{Proposition}

\theoremstyle{remark}
\newtheorem{remark}[theorem]{Remark}

\theoremstyle{definition}

\numberwithin{equation}{section}

\usetikzlibrary{datavisualization}
\usetikzlibrary{datavisualization.formats.functions}
\usetikzlibrary{decorations.pathmorphing}
\usetikzlibrary{decorations.pathreplacing}
\usetikzlibrary{matrix, arrows, calc}
\usetikzlibrary{arrows,decorations.pathmorphing,backgrounds,positioning,fit}

\title[A short FE code for 2D FL]{ A short FE implementation for a 2d homogeneous Dirichlet problem of a Fractional Laplacian
}
\author[G. Acosta, F. M. Bersetche and J. P. Borthagaray]{Gabriel Acosta, 
Francisco M. Bersetche and Juan Pablo Borthagaray}

\address[G. Acosta, F. M. Bersetche and J. P. Borthagaray]{IMAS - CONICET and 
Departamento de Matem\'a\-tica, FCEyN - Universidad de Buenos Aires, Ciudad 
Universitaria, Pabell\'on I  (1428) Buenos Aires, Argentina.}

\email[G. Acosta]{gacosta@dm.uba.ar}

\urladdr[G. Acosta]{http://mate.dm.uba.ar/~gacosta/}

\email[F. M. Bersetche]{fbersetche@dm.uba.ar}

\email[J.P. Borthagaray]{jpbortha@dm.uba.ar}

\thanks{This work has been partially supported by CONICET, ANPCYT and UBA under 
grants PIP 2014 1220130100034CO, PICT 2014-1771 and UBACYT 20020130100205BA
}
\subjclass[2010]{65N30, 35R11}
\keywords{Finite Elements, Fractional Laplacian, Nonlocal Operators}

\begin{document}
\begin{abstract}
In \cite{AcostaBorthagaray}, a complete $n$-dimensional
finite element analysis of the homogeneous Dirichlet problem associated 
to a fractional Laplacian was presented.   Here we provide a comprehensive and 
simple 2D {\it MATLAB}\textsuperscript{\textregistered} finite element code
for such a problem. The code is accompanied with a basic discussion of the theory 
relevant in the context. The main program is written in about 80 lines  
and can be easily modified to deal with other kernels 
as well as with time dependent problems. The 
present work fills a gap by providing an input for a large number of 
mathematicians and scientists interested in numerical approximations
of solutions of a large variety of problems involving nonlocal phenomena in two-dimensional space.
 \end{abstract}
\maketitle

\section{Introduction}
The Finite Element Method (FEM) is one of the preferred numerical tools
in scientific and engineering communities. It counts with a solid and 
long established theoretical foundation, mainly in the linear case
of second order \emph{elliptic}  partial differential
equations. These kind of operators, with the Laplacian as a canonical example,
are involved in modeling \emph{local} diffusive 
processes. On the other hand, nonlocal or \emph{anomalous} diffusion models
have increasingly impacted upon a number of important
areas in science. Indeed, non-local
formulations can be found in
physical and social contexts, modeling as diverse phenomena as human locomotion
in relation to crime diffusion \cite{crime}, electrodiffusion of ions within
nerve cells \cite{Langlands} or machine learning \cite{Rosasco}.

The Fractional Laplacian (FL)  is among the most 
prominent examples of a non-local operator. For $0<s<1$, it is defined
as 
\begin{equation}
(-\Delta)^s u (x) = C(n,s) \mbox{ p.v.} \int_\rn \frac{u(x)-u(y)}{|x-y|^{n+2s}} \, dy,
\label{eq:op_fraccionario}
\end{equation}
where
$$C(n,s) = \frac{2^{2s} s \Gamma(s+\frac{n}{2})}{\pi^{n/2} \Gamma(1-s)} $$
is a normalization constant. The FL, given by \eqref{eq:op_fraccionario}, is one of the 
simplest  pseudo-differential operators and can also be 
regarded as the infinitesimal generator of a $2s$-stable L\'evy process 
\cite{bertoin}. 

Given a function $f$ defined in a bounded domain 
$\W$, the homogeneous Dirichlet problem  associated to the FL reads: find $u$ such that
\begin{equation}
\left\lbrace
  \begin{array}{rl}
      (-\Delta)^s u = f & \mbox{ in }\W, \\
      u = 0 & \mbox{ in }\W^c . \\
      \end{array}
    \right.
\label{eq:fraccionario}
\end{equation}

In contrast to elliptic PDEs, numerical developments for problems involving  this
non-local operator, even in simplified contexts, are seldom found in the literature. 
The reason for that is related to two major challenging
tasks usually involved in its numerical treatment:  the
handling of highly singular kernels and the need to cope with an unbounded region
of integration.  This is precisely  the case of \eqref{eq:fraccionario},
for which just a few numerical methods have been proposed. 
Effectively implemented \emph{in one space dimension}, we mention, 
for instance:
a finite difference scheme by Huang and Oberman \cite{HuangOberman},
 a FE approach  developed by D'Elia and Gunzburger \cite{DEliaGunzburger} 
that relies on a volume-constrained version of the non-local operator and 
a simple one-dimensional  spectral approach \cite{AcostaBorthagarayBrunoMaas}. 
We refer the reader to \cite{AcostaBorthagaray} for a more detailed account of these
schemes and a discussion on other fractional diffusion operators on bounded domains and their 
discretizations.

To the best of the authors' knowledge,
numerical computations for \eqref{eq:fraccionario} in higher dimensions have become available 
only recently 
 \cite{AcostaBorthagaray}. In  that paper a complete $n$-dimensional
 finite  element analysis for the FL  has been carried out, including
 regularity of solutions of 
\eqref{eq:fraccionario} in standard and weighted fractional spaces.  
Moreover, the convergence for piecewise linear elements 
 is proved with optimal order for both uniform and graded meshes. 

In that work there are presented error bounds in the energy norm and numerical experiments (in 2D), demonstrating an accuracy of the order of $h^{1/2} \log h$
and $h \log h$ for solutions obtained by means of uniform and graded meshes, 
respectively.

The present  article can be seen as a complementary work to 
\cite{AcostaBorthagaray}, providing a short and simple {\it MATLAB}\textsuperscript{\textregistered} FE 
code coping with the homogeneous Dirichlet problem \eqref{eq:fraccionario}. 

In  \cite{50lines} a 
{\it MATLAB}\textsuperscript{\textregistered} implementation for \emph{linear} finite elements 
and \emph{local} elliptic operators is presented in a concise way.  We tried 
to emulate as much as possible that spirit in the non-local context. 
Notwithstanding that and in spite of our efforts, some intrinsic technicalities 
make our code inevitably slightly longer and more complex than that. Just to give 
a hint about it, we take a glimpse in advance  at
the nonlocal stiffness matrix $K$. It involves expressions of the type
\begin{equation}
 \label{eq:ejemplitototal}
\int_{\R^2}\int_{\R^2} \frac{(\phii_i(x) - 
\phii_i(y))(\phii_j(x)-\phii_j(y))}{|x-y|^{2+2s}} \, dx dy, 
\end{equation}
where $\phii_i,\phii_j$ are arbitrary nodal basis functions associated
to a triangulation $\mathcal{T}$. Two 
difficulties become 
apparent in the calculation of \eqref{eq:ejemplitototal}. First, 
at the element level, computing \eqref{eq:ejemplitototal} leads 
to terms like 
\begin{equation}
 \label{eq:ejemplito}
\int_{T}\int_{\tilde{T}} \frac{(\phii_i(x) - 
\phii_i(y))(\phii_j(x)-\phii_j(y))}{|x-y|^{2+2s}} \, dx dy, 
\end{equation}
for arbitrary  pairs  $T,\tilde{T}\in \mathcal{T}$. 
If $T$ and $\tilde{T}$ are not neighboring then the integrand
in \eqref{eq:ejemplito} is a regular function and can be integrated
numerically in a standard fashion. On the other hand, if
$T\cap \tilde{T}\neq \emptyset$ an accurate 
algorithm to compute  \eqref{eq:ejemplito} is not easy to devise. 
Fortunately, \eqref{eq:ejemplito} bears some resemblances to 
typical integrals appearing in the  Boundary Element 
Method \cite{SauterSchwab} and we extensively exploit this fact. 
Indeed, a basic and well known technique in the BEM community
is to rely on Duffy-type transforms. This approach leads us to the decomposition of such integrals into two parts:  a highly singular but 
explicitly integrable part and a smooth, numerically treatable part.
We use this method to show how \eqref{eq:ejemplito} 
can be handled with an arbitrary degree of precision (this is carefully treated 
in Appendices 
\ref{sec:nontouching}, \ref{sec:vertex}, \ref{sec:edge}, \ref{sec:identical}). 

Yet another difficulty is hidden in the calculation of $K$.  Although $\Omega$ is a bounded domain and the number of potential unknowns is always finite, \eqref{eq:ejemplitototal} involves a computation in $\R^2 \times \R^2$.  In particular, in the homogeneous setting, 
we need to accurately compute the function 
\begin{equation}
 \label{eq:ejemplito2}
\int_{\Omega^c} \frac{1}{|x-y|^{2+2s}} \,  dy, 
\end{equation}
for any $x\in \Omega$. That, of course,  can be hard to achieve 
for a domain with a complex boundary. Nonetheless, introducing an
extended secondary mesh, as it is explained in Section \ref{sec:preliminar}, 
it is possible to reduce such problem to a simple case in which $\partial \Omega$ 
is a circle. We show that in this circumstance a computation of
\eqref{eq:ejemplito2} can be both fast and accurately delivered (see 
also Appendix \ref{sec:comp}). Remarkably,  this simple idea applies in 
arbitrary space dimensions.

Regarding the code itself,  our main 
concern has been to keep a compromise between readability and 
efficiency. 
First versions of our code were plainly readable
but too slow to be satisfactory. In the code offered here many
computations have been vectorized and a substantial speed up 
gained, sometimes at the price of losing (hopefully
not too much) readability.

Last but not least, the full program is available from the authors upon request, so that 
the reader can avoid retyping it. Small modifications of the base code may make it  
usable for  dealing with many different problems. It has been  successfully used 
in several contexts such as eigenvalue computations and time dependent problems
(considering semi and full fractional settings), among others.

The paper is organized as follows. In Section 
\ref{sec:prelim}, we review appropriate fractional spaces and regularity results for \eqref{eq:fraccionario}. 
Section \ref{sec:preliminar} deals with basic aspects of the
FE setting. The data structure is carefully discussed in Section \ref{sec:codigo1}
and the main loop of the code is described in Section \ref{sec:mainloop}. 
Section \ref{sec:numex}, in turn, shows a numerical example for which a nontrivial (i.e. with a non constant source 
term $f$) solution is explicitly known. Moreover, the e.o.c. in $L^2(\W)$ is presented for some values of $s$.
These numerical results are in very good agreement with those expected
by using standard duality arguments together with the theory given in 
\cite{AcostaBorthagaray}. Appendix \ref{sec:quadrature} may be found rather technical for people not coming
from the Boundary Element community and deals with  the quadrature rules used in 
each singular case.  Appendices \ref{sec:fun} and \ref{sec:data} describe respectively 
auxiliary functions and data 
used along the program. Finally, the full code, including the line numbers, is 
exhibited in Appendix \ref{sec:maincode}.

\section{Function spaces and regularity of solutions} 
\label{sec:prelim}
Given an open set $\W \subset \rn$  and $s \in(0,1)$, define the fractional 
Sobolev space $H^s(\W)$ as
	\[
		H^s(\W) = \left\{ v \in L^2(\W) \colon |v|_{H^s(\W)} < \infty 
\right\},
	\]
	where $|\cdot|_{H^s(\W)}$ is the Aronszajn-Slobodeckij seminorm
	\[
		|v|_{H^s(\W)}^2 = \iint_{\W^2} 
		\frac{|v(x)-v(y)|^2}{|x-y|^{n+2s}} \, dx \, dy. 
	\]
	It is evident that $H^s(\W)$ is a Hilbert space endowed with the norm 
$\|\cdot\|_{H^s(\W)} = \|\cdot\|_{L^2(\W)} + |\cdot|_{H^s(\W)} .$
	Moreover, consider the bilinear form $\langle \cdot, \cdot \rangle 
$ on $H^s(\W),$ 
	\begin{equation}
		%\langle \cdot, \cdot \rangle \colon H^s(\W) \times H^s(\W) \to 
%\R,	\quad 
		\langle u, v \rangle_{H^s(\W)} = \iint_{\W^2} 
\frac{(u(x)-u(y))(v(x)-v(y))}{|x-y|^{n+2s}} \, dx \, dy.
	\label{eq:prod_int}
	\end{equation}
Let us also define the space of functions supported in $\W$, 
\[
	\widetilde{H}^s (\W) = \left\{ v \in H^s(\rn) \colon \text{ supp } v 
\subset \bar{\W} \right\}.
\]	
This space may be defined through interpolation,
\[ 
	\widetilde{H}^s (\W) = \left[L^2(\W), H^1_0(\W) \right]_{s}.
\]
Moreover, depending on the value of $s$, different characterizations of this 
space are available.
If $s<\frac12$ then $\widetilde{H}^s (\W)$ coincides with $H^s (\W)$, and 
if $s>\frac12$ it may 
be characterized as the closure of $C^\infty_0(\W)$ with respect to the 
$|\cdot|_{H^s (\W)}$ norm. In the latter case, it is also customary to denote it 
by $H^s_0(\W)$. The particular case of $s=\frac12$ gives raise to the 
Lions-Magenes space $H^{\frac12}_{00}(\W)$, 
which can be characterized by
\[
	H^{\frac12}_{00} (\W) = \left\{ v \in H^{\frac12}(\W) \colon \int_\W 
\frac{v(x)^2}{\text{dist}(x,\partial \W)} \, dx < 
	\infty \right\}.
\]
Note that the inclusion $H^{\frac12}_{00}(\W) \subset H^{\frac12}_{0}(\W) 
= H^{\frac12}(\W)$ is strict. We also need to introduce 
the dual space of $\widetilde{H}^s(\W)$,  denoted with the standard
negative exponent $H^{-s}(\Omega)$.

It is apparent that the form $\langle \cdot, \cdot \rangle_{H^s(\rn)}$ (recall 
\eqref{eq:prod_int})
induces a norm on $\widetilde{H}^s(\W),$ because of the following well 
known result.
	\begin{proposition}[Poincar\'e inequality]
		There is a constant $c=c(\W,n,s)$ such that 
		\begin{equation*} \label{eq:poincare}
			\| v \|_{L^2(\W)} \leq c |v|_{H^s(\rn)}  \quad \forall v 
\in \widetilde{H}^s(\W) .
		\end{equation*}
	\end{proposition}

Finally, Sobolev spaces of order grater than 1 are defined in the following way: given $k 
\in \mathbb{N}$, then
$$H^{k+s} (\W) = \left\{ v \in H^k(\W) \colon |D^\alpha v | \in H^s(\W) \, 
\forall \alpha \mbox{ with } |\alpha| = k \right\} ,$$
furnished with the norm
$$ \| v \|_{H^{k+s} (\W)} = \| v \|_{H^k(\W)} + \sum_{|\alpha| = k } | D^\alpha 
v |_{H^s(\W)}.$$

Weak solutions of \eqref{eq:fraccionario} are straightforwardly defined 
multiplying by a test function and integrating by parts. Indeed, the weak 
formulation of \eqref{eq:fraccionario} reads:  find $u \in \widetilde{H}^s(\W)$ 
such that
\begin{equation}
\frac{C(n,s)}{2}\langle u, v \rangle_{H^s(\rn)} = \int_\W f v, \quad v \in 
\widetilde{H}^s(\W) .
\label{eq:debil}
\end{equation}
Notice that the inner product  
\begin{equation}
 \langle u, v \rangle_{H^s(\rn)} =\iint_{\R^n\times \R^n} 
\frac{(u(x)-u(y))(v(x)-v(y))}{|x-y|^{n+2s}} \, dx \, dy.
	\label{eq:prod_intGLOBAL}
\end{equation}
involves integrals in $\R^n$.

From now on, we assume $f \in H^r(\W)$ for some $r \ge -s$. Existence 
and uniqueness of solutions in $\tilde{H}^s(\W)$ and well-posedness of problem 
\eqref{eq:debil} are immediate consequences of the Lax-Milgram lemma. Moreover, 
the following regularity result is valid \cite{Grubb, VishikEskin}:
\begin{theorem}
\label{teo:grubb-Eskin}
Let $u \in \widetilde H^s(\W)$ be the solution to \eqref{eq:debil}. If $\partial \Omega$
is of $C^\infty$ class, then 
\[ u \in
\begin{cases}
H^{2s + r}(\W) & \text{if } s + r < 1/2, \\
H^{s + 1/2 - \eps}(\W) \ \forall \eps > 0 & \text{if } s + r \ge 1/2 .
\end{cases}
\]
\end{theorem}

\begin{remark} The previous theorem implies that, 
independently of the  regularity of the right hand side 
function $f$, solutions should not be expected to have 
derivatives of order greater than $s+1/2$  in $L^2(\W)$. This is a 
consequence of the behavior of solutions near the boundary of $\W$: the quotient 
$u(x)/d(x,\partial\W)^s$ can be shown to be finite for $x\sim \partial \W$ (see, 
for example \cite{RosOtonSerra}). Knowledge of this singularity was exploited in 
\cite{AcostaBorthagaray}, where problem \eqref{eq:debil} was set up in the 
framework of weighted Sobolev spaces and solutions were proved to have 
$1+s-\eps$ derivatives in a suitable space if the right hand side function 
belongs to $C^{1-s}(\W)$. See that work for further details. 
\end{remark}

%%%%%%%%%%%%%%%%%%%%%%%%%%%%%%%%%%%%%%%%%%%%%%%%%%%%%%%%%%%%%%%%%%%%%%%%%%%%%%%%
%%%%%

\section{FE setting}
\label{sec:preliminar}

Consider  an admissible triangulation $\mathcal{T}$ of $\W$ consisting of 
$N_\mathcal{T}$ elements. For the discrete space $\V_h$, we take
 standard \emph{continuous} piecewise linear elements over $\mathcal{T}.$ 
With the usual notation, we introduce the nodal basis 
$\{ \phii_1, \ldots, \phii_N\}\subset \V_h$ corresponding to the internal nodes $\{ x_1, \ldots, x_N\}$, 
that is $\phii_i(x_j)=\delta_i^j$. Given an element 
$T \in \mathcal{T}$, we denote by $h_T$ and $\rho_T$ its 
diameter and inner radius, respectively. As customary, we write $h=\max_{T\in \mathcal{T}} h_T.$
The family of triangulations considered is assumed to be shape-regular, 
namely,  there exists $\sigma>0$ independent of $\mathcal{T}$ such that
\[
h_T \le \sigma \rho_T \text{ for all } T \in \mathcal{T}.
\]
In this context, the  discrete analogous of \eqref{eq:prod_intGLOBAL}
reads: 
find $u_h \in \V_h$ such that 
\begin{equation}
\frac{C(n,s)}{2}\langle u_h, v_h \rangle_{H^s(\rn)} = \int_\W f v_h, \quad v_h 
\in \V_h ,
\label{eq:debil_discreto}
\end{equation}
providing a \emph{conforming}\footnote{Notice that even $P_0$ elements are conforming for $0<s<1/2$. We restrict ourselves to
continuous $P_1$ in order to give an unified conforming approach for any $0<s<1$.} FEM for any $0<s<1$.  

Writing the 
discrete solution as $u_h = \sum_j u_j \phii_j$, problem 
\eqref{eq:debil_discreto} is equivalent to solving the linear system
\begin{equation}
KU = F,
\label{eq:sistema}
\end{equation} 
where the coefficient matrix $K=(K_{ij}) \in \R^{N\times N}$ and the right-hand 
side $F=(f_j)\in \R^N$ are defined by
\[
K_{ij} = \frac{C(n,s)}{2} \langle \phii_i, \phii_j \rangle_{H^s(\rn)}, \quad 
f_j= \int_\W f \phii_j,
\]
and the unknown is $U = (u_j) \in \R^N$.

\begin{figure}[ht]
$\begin{array}{cc}
	\includegraphics[width=0.35\textwidth]{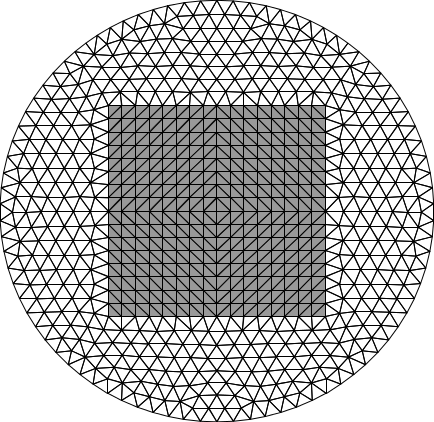}
\end{array}$
		\caption{A square domain $\W$ (gray) and an auxiliary ball containing 
		it. Regular triangulations $\mathcal{T}$ and $\mathcal{T}_A$ for 
		$\W$ and $B\setminus \Omega$ are shown. The final symmetry
		of the admissible triangulation $\tilde{\mathcal{T}}=\mathcal{T} \cup \mathcal{T}_A$, exhibited in the example, is not relevant.}
\label{fig:integrafuera}
\end{figure}

The \emph{fractional} stiffness matrix $K$ is symmetric and positive definite, so that 
\eqref{eq:sistema} has a unique solution. 
Notice that the integrals in the inner product involved in 
computation of $K_{ij}$ should be carried over $\R^n$. For this reason we
find it useful to consider a ball
$B$ containing $\W$ and such that the distance from  $\bar \Omega$ to $B^c$ is an 
arbitrary positive number.  As it is explained in Appendix  \ref{sec:comp}, this is needed in order to avoid 
difficulties caused by lack of symmetry when dealing with the integral 
over $\W^c$ when $\W$ is not a ball. Together with $B$,  we introduce an auxiliary 
triangulation ${ \mathcal{T}_A}$ on $B\setminus\W$ such that the 
complete triangulation $\tilde{\mathcal{T}}$ over $B$ (that is $\tilde{\mathcal{T}}=
\mathcal{T}\cup \mathcal{T}_A$) is admissible (see Figure \ref{fig:integrafuera}). 

Let us call $ N_{ \mathcal{\tilde T}}$ the number of elements on the 
triangulation of $B$. Then, defining for $1\le \ell, m 
\le N_{\mathcal{\tilde T}}$ and $1\le \ell \le 
N_{\mathcal{\tilde T}}$ 
\begin{equation}
\begin{split}
I^{i,j}_{\ell,m} & = \int_{T_\ell}\int_{T_m} \frac{(\phii_i(x) - 
\phii_i(y))(\phii_j(x)-\phii_j(y))}{|x-y|^{2+2s}} \, dx dy, \\
J^{i,j}_\ell & = \int_{T_\ell}\int_{B^c} 
\frac{\phii_i(x)\phii_j(x)}{|x-y|^{2+2s}} \, dy dx,  \label{eq:integrales}
\end{split}
\end{equation}
we may write
\[ 
K_{ij} =  \frac{C(n,s)}{2} \sum_{\ell = 1}^{N_{\mathcal{\tilde T}}} \left( 
\sum_{m=1}^{N_{\mathcal{\tilde T}}} I^{i,j}_{\ell,m} + 2 J^{i,j}_\ell \right).
\]
As mentioned above, the computation of each integral $I_{\ell, m}^{i,j}$ and 
$J_\ell^{i,j}$ is  challenging for different reasons: the former involves a singular integrand
 if 
$\overline{T_\ell}\cap\overline{T_m}\neq \emptyset$ 
(Appendices 
\ref{sec:vertex}, \ref{sec:edge}, \ref{sec:identical} are devoted to handle it) while the latter needs to 
be calculated on an unbounded domain. In this case notice that
$$
J^{i,j}_\ell= \int_{T_\ell}
\phii_i(x)\phii_j(x)
\psi(x) \, dx,
$$
with $\psi(x):= \int_{B^c} \frac{1}{|x-y|^{2+2s}} \, dy$.
Therefore all we need is an accurate computation of 
$h(x)$ for each quadrature point used
in $T_\ell\subset \bar{\Omega}$ (notice that $h(x)$ is a smooth function up to the
boundary of $\Omega$ since $|x-y|>dist(\bar \Omega, B^c)>0$). 

Taking this into account, we observe that it is possible to take advantage of the fact that
$h(x)$ is a radial function that can be either quickly computed
on the fly or even precomputed with an arbitrary degree of precision (see Appendix  \ref{sec:comp} for a full treatment of $h(x)$).

For the reader's convenience we finish this section with Table \ref{tablaTotal}, containing
some handy notations. 
\begin{table}[h]
\centering
\caption{Main Variables}
\label{tablaTotal}
\begin{tabular}{|c|c|}
\hline
Notation            & Meaning   \\ \hline
$\mathcal{T},\mathcal{T}_A,\tilde{\mathcal{T}}$  & Meshes: of $\Omega$, $B\setminus \Omega$ and $B$ resp.\\ \hline
$\mathcal{N}$    & Nodes of $\mathcal{T}$  \\ \hline   
$\mathcal{E}$   &   Edges of $\mathcal{T}$  \\ \hline
$\mathcal{B}$   &   Boundary edges  of $\mathcal{T}$   \\ \hline
$N_{\mathcal{T}}$ & $\#\mathcal{T}$\\ \hline
$N_{\mathcal{N}}$ & $\#\mathcal{N}$\\ \hline
$N_{\mathcal{B}}$ & $\#\mathcal{B}$\\ \hline
\end{tabular}
\end{table}

%%%%%%%%%%%%%%%%%%%%%%%%%%%%%%%%%%%%%%%%%%%%%%%%%%%%%%%%%%%%%%%%%%%%%%%%%%%%%%%%
%

%%%%%%%%%%%%%%%%%%%%%%%%%%%%%%%%%%%%%%%%%%%%%%%%%%%%%%%%%%%%%%%%%
\begin{section}{Data structure and auxiliary variables}
\label{sec:codigo1}
%%%%%%%%%%%%%%%%%%%%%%%%%%%%%%%%%%%%%%%%%%%%%%%%%%%%%%%%%%%%%%%%%

We assume that the mesh $\mathcal{T}$ has been generated in advance\footnote{For the 
sake of convenience an stored example mesh -as well as  
a suitable mesh generator- is provided together with the source code.}. The information related to $\mathcal{T}$ should
be encoded in some specific variables \verb+p+, \verb+t+, 
\verb+bdrynodes+, \verb+nt_aux+ \verb+nf+ \verb+R+, as follows: 

\begin{itemize}
\item
\verb+p+ is a $2 \times N_{\mathcal{N}}$ array, such that \verb+p(:,n)+ are the 
coordinates of the \verb+n+-th node.
\item
 \verb+t+ is a $N_{\tilde{\mathcal{T}}} \times 3$ index array, and 
 \verb+t(l,:)+ are the 
indices of the vertices of $T_l$. Triangles belonging to $\mathcal{T}_{A}$ must be listed at the end. 
\item 
  \verb+nt_aux+ $= \#\mathcal{T}_{A}$.   
\item
  \verb+bdrynodes+ is an index column vector listing the 
  nodes lying on  $\partial \Omega$. 
\item
  \verb+nf+ is an index column vector contiaining the free nodes 
  (those in $\Omega$). 
\item
  \verb+R+ the radius of $B$.

\end{itemize}
These data have to be available in the {\it MATLAB}\textsuperscript{\textregistered} workspace before
 the execution of the main code. 
 
 Next, we begin by creating some variables that refer to problem 
\eqref{eq:fraccionario}:

\begin{Verbatim}[fontsize=\small]
s = 0.5;
f = @(x,y) 1; 
cns = s*2^(-1+2*s)*gamma(1+s)/(pi*gamma(1-s));
load(`data.mat');
\end{Verbatim}

Here, \verb+s+ is the order of the fractional Laplacian involved, \verb+f+ is a 
function handle containing the volume force (which as an example we have set to be $f \equiv 1$), and \verb+cns+ is equal to the 
constant $C(n,s)$ previously defined. 

In order to compute the stiffness matrix we need to estimate the 
bilinear form $\langle \cdot, \cdot \rangle_{H^s(\rn)}$ evaluated at the nodal 
basis through an appropriate quadrature rule. 

To perform an efficient vectorized computation, we require some pre-calculated 
data, given in the file \verb+data.mat+. 
This file contains information about nodes and weights for the quadratures 
performed throughout the code.
The content of \verb+data.mat+ is listed in Table 
\ref{tablaData} and further details can be found in 
Appendix \ref{sec:data}.

\begin{table}[h]
\centering
\caption{Variables stored in \texttt{data.mat}}
\label{tablaData}
\begin{tabular}{|l|l|l|l|}
\hline
Name            & Size  &\begin{tabular}[c]{@{}l@{}}Used as input \\ in 
function: \end{tabular}                                          & Description   
                                                                        \\ 
\hline
\verb+p_cube+  & 27x3  & \begin{tabular}[c]{@{}l@{}}\verb+vertex_quad+\\ 
\verb+edge_quad+\end{tabular} & \begin{tabular}[c]{@{}l@{}}Quadrature points \\ 
over $[0,1]^3$\end{tabular}           \\ \hline
\verb+p_T_6+       & 6x2   & \begin{tabular}[c]{@{}l@{}}None (used in non-\\ 
touching case)\end{tabular}     & \begin{tabular}[c]{@{}l@{}}Quadrature points 
\\ over $\hat{T}$\end{tabular}           \\ \hline
\verb+p_T_12+ & 12x2  &       \verb+comp_quad+                         & 
\begin{tabular}[c]{@{}l@{}}Quadrature points\\  over $\hat{T}$\end{tabular}      
     \\ \hline
\verb+p_I+       & 9x1   & \begin{tabular}[c]{@{}l@{}}\verb+comp_quad+\\ 
\verb+triangle_quad+\end{tabular}                                                
                     & \begin{tabular}[c]{@{}l@{}}Quadrature points\\ over 
$[0,1]$\end{tabular}              \\ \hline
\verb+w_I+       & 9x1   & \verb+comp_quad+                                      
                         & \begin{tabular}[c]{@{}l@{}}Quadrature weights \\ 
associated to \verb+p_I+\end{tabular} \\ \hline
\verb+phiA+     &      &                                                         
                     &                                                           
          \\
\verb+phiB+     & 9x36 &  \begin{tabular}[c]{@{}l@{}}None (used in non-\\ 
touching case)\end{tabular}  & See Appendix 
\ref{sec:nontouchingdata}                        \\
\verb+phiD+     &       &                                                        
                         &                                                       
                                \\ \hline
\verb+vpsi1+    & 25x27 & \verb+vertex_quad+                                     
                        & See Appendix 
\ref{sec:vertexdata}                                                           
           \\
\verb+vpsi2+    &       &                                                        
                         &                                                       
                                \\ \hline
\verb+epsi1+    &      &                                                         
        &                                                                       
\\
\verb+epsi2+    &       &                                                        
                         &                                                       
                                \\
\verb+epsi3+    & 16x27 &         \verb+edge_quad+                               
                                         &   See 
Appendix \ref{sec:edgedata}                            
                                       \\
\verb+epsi4+    &       &                                                        
                         &                                                       
                                \\
\verb+epsi5+    &       &                                                        
                         &                                                       
                                \\ \hline
\verb+tpsi1+    &       &                                                        
       &                                                                         
                                                                     \\
\verb+tpsi2+    &9x9&   \verb+triangle_quad+                                     
               &                  See Appendix 
\ref{sec:identicaldata}                                                        
                 \\
\verb+tpsi3+    &       &                                                        
                         &                                                       
                                \\ \hline
\verb+cphi+    & 9x12 & \verb+comp_quad+                                         
                       & See Appendix 
\ref{sec:compdata}                                                             
         \\ \hline
\end{tabular}
\end{table}

As mentioned before, some auxiliary elements are added to the original mesh 
in order to have a triangulation on a ball $B$ containing $\W$ (see 
Figure \ref{fig:integrafuera}). The 
nodes in this auxiliary domain $B \setminus \W$ are regarded as Dirichlet nodes.

Next, we define some mesh parameters and set to zero the factors involved in 
equation \eqref{eq:sistema}. The following lines do not need extra explanation 
beyond the in-line comments:
\begin{Verbatim}[fontsize=\small]
nn = size(p,2); % number of nodes
nt = size(t,1) % number of elements
uh = zeros(nn,1); % discrete solution
K  = zeros(nn,nn); % stiffness matrix
b  = zeros(nn,1); % right hand side
\end{Verbatim}

Then, the measures of all the elements in the mesh are calculated: 
\begin{Verbatim}[fontsize=\small]
area = zeros(nt,1);
for i=1:nt
    aux = p( : , t(i,:) );
    area(i) = 0.5.*abs(...
	 det([ aux(:,1) - aux(:,3)  aux(:,2) - aux(:,3)]) );
end
\end{Verbatim}
So, \verb+area+ is a vector of length $N_{\mathcal{\tilde T}}$ satisfying 
$\verb+area(l)+ = |T_{\verb+l+}|$, $\verb+l+ \in \{1, ... , N_{\mathcal{\tilde 
T}}\}$.  

The quadratures we employ to compute the integrals $I_{\ell,m}^{i,j}$ (defined 
in \eqref{eq:integrales}) depend on whether the elements $T_{\ell}$ and $T_m$ 
coincide or their intersection is an edge, a vertex or empty.  Therefore, it is 
important to distinguish theses cases in an efficient way. We construct a data 
structure called \verb+patches+ as follows, using a linear number of operations:
\begin{Verbatim}[fontsize=\small]
deg = zeros(nn,1);
for i=1:nt
    deg( t(i,:) ) = deg( t(i,:) ) + 1;
end
patches = cell(nn , 1);
for i=1:nn
    patches{i} = zeros( 1 , deg(i) );
end
for i=1:nt
    patches{ t(i,1) }(end - deg( t(i,1) ) + 1) = i;
    patches{ t(i,2) }(end - deg( t(i,2) ) + 1) = i;
    patches{ t(i,3) }(end - deg( t(i,3) ) + 1) = i;
    deg( t(i,:) ) = deg( t(i,:) ) - 1;
end
\end{Verbatim}
The output of this code block is a $N_{\tilde{\mathcal N}} \times 1$ cell, 
called  \verb+patches+, such that \verb+patches{n}+ is a vector containing the 
indices of all the elements in the neighborhood of the node \verb+n+. 

\end{section}

%%%%%%%%%%%%%%%%%%%%%%%%%%%%%%%%%%%%%%%%%%%%%%%%%%%%%%%%%%%%

\begin{section}{Main loop}
\label{sec:mainloop}
One of the main challenges to build up a FE implementation to problem 
\eqref{eq:fraccionario} is to assemble the stiffness matrix in an efficient 
mode. 
Independently of whether the supports of two given basis functions $\phii_i$ and 
$\phii_j$ are disjoint, the interaction $\langle \phii_i, \phii_j 
\rangle_{H^s(\rn)}$ is not null. This yields a paramount difference between FE 
implementations for the classical and the fractional Laplace operators; in the 
former the stiffness matrix is sparse, while in the latter it is full.
Therefore, unless some care is taken,
%for an straightforward implementation, 
the amount of memory required and the number of operations needed to assembly 
the stiffness matrix increases quadratically with the number of nodes.
Due to this, the code we present takes advantage of vectorized operations as 
much as possible.   

Moreover, as the computation of the entries of the stiffness matrix requires 
calculating integrals on \emph{pairs} of elements, it is required to perform a 
double loop. It is simple to check  the identity $I_{\ell,m}^{i,j} = I_{m, 
\ell}^{i,j}$ for all $i, j, \ell, m$, and therefore it is enough to carry the 
computations only for the pairs of elements $T_\ell$ and $T_m$ with $\ell \le 
m$.

 In the following lines we preallocate memory and create the auxiliary index 
array \verb+aux_ind+ (to be used in code line 58).  
\begin{Verbatim}[fontsize=\small]
vl = zeros(6,2); 
vm = zeros(6*nt,2); 
norms = zeros(36,nt);
ML = zeros(6,6,nt);
empty = zeros(nt,1);   
aux_ind = reshape( repmat( 1:3:3*nt , 6 , 1 ) , [] , 1 );
empty_vtx = zeros(2,3*nt);
BBm = zeros(2,2*nt);
\end{Verbatim}

The main loop goes through all the elements $T_{\ell}$ of the mesh of $\W$, 
namely, $1 \le \ell \le N_\mathcal{T}$.  Observe that auxiliary elements are 
excluded from it. 
Fixed $\ell$, the first task is to classify all the mesh elements $T_m$ ($1 \le m \le 
N_{\tilde{\mathcal{T}}}$, $m \ne \ell$) according to whether 
$\overline{T_{\ell}} \cap \overline{T_m }$ is empty, a vertex or an edge. This 
is accomplished employing a linear number of operations by using the 
\verb+patches+ data structure as follows: 
\begin{Verbatim}[fontsize=\small]
edge =  [ patches{t(l,1)} patches{t(l,2)} patches{t(l,3)} ];
[nonempty M N] = unique( edge , 'first' );
edge(M) = [];
vertex = setdiff(  nonempty , edge  );          
ll = nt - l + 1 - sum( nonempty>=l );
edge( edge<=l ) = []; 
vertex( vertex<=l ) = []; 
empty( 1:ll ) = setdiff_( l:nt , nonempty ); 
empty_vtx(: , 1:3*ll) = p( : , t( empty(1:ll) , : )' );
\end{Verbatim}

At this point,  \verb+ll+ is the number of elements --including the auxiliary 
ones-- whose intersection with $T_{\ell}$ is empty and have not been visited yet 
(namely, those with index \verb+m+$>$\verb+l+). 
By considering only the elements with index greater than $\ell$, we are taking 
advantage of the symmetry of the stiffness matrix. The arrays
\verb+empty+, \verb+vertex+ and  \verb+edge+ contain the indices of all those elements whose 
intersection with $T_{\ell}$ is empty, a vertex or an edge respectively, and 
have not been computed yet. In \verb+empty_vtx+ we store the coordinates of the 
vertices of the triangles indexed in \verb+empty+.

Then, the code proceeds to assemble the right hand side vector in equation 
\eqref{eq:sistema} 
\begin{Verbatim}[fontsize=\small]
nodl = t(l,:);
xl = p(1 , nodl); yl = p(2 , nodl);
Bl = [xl(2)-xl(1) yl(2)-yl(1); xl(3)-xl(2) yl(3)-yl(2)]'; 
b(nodl) = b(nodl) + fquad(area(l),xl,yl,f); 
\end{Verbatim}
Here, \verb+nodl+ stores the indices of the vertices of $T_{\ell}$; \verb+xl+ 
and \verb+yl+ are the $x$ and $y$ coordinates of these vertices, respectively. 
The element $T_\ell$ is the image of a reference element $\hat T$ via an affine 
transformation, 
\[
( \hat x , \hat y ) \mapsto \verb+Bl+ ( \hat x , \hat y) + (\verb+xl(1)+ , 
\verb+yl(1)+) .
\]

Recall that \verb+b+ stores the numerical approximation to the right hand side 
vector from equation \eqref{eq:sistema}, namely, $\verb+b(j)+\approx \int_\W f 
\phii_j.$
The routine \verb+fquad+ uses a standard quadrature rule, interpolating $f$ on 
the edge midpoints of $T_l$ (see Appendix \ref{sec:fun}).

\begin{remark} \label{rem:numeracion}
Let $1 \le \ell, m \le \mathcal{N}_{\tilde T}$. When computing 
$I_{\ell,m}^{i,j}$ or $J_{\ell}^{i,j}$, 
the basis function indices $i$ and $j$ do not refer to a global numbering but 
to a local one. This means, for example, that if $\overline{T_\ell} \cap 
\overline{T_m} = \emptyset$, then $1 \le i, j \le 6$. See Remark 
\ref{rem:numeracion_ap} for details on this convention. 
\end{remark}

\subsection{Identical elements}
The first interaction to be computed by the code corresponds to the case 
$m=\ell$ in \eqref{eq:integrales}. The values calculated are assembled in the 
stiffness matrix \verb+K+.
\begin{Verbatim}[fontsize=\small]
K(nodl, nodl) = K(nodl, nodl) +...
triangle_quad(Bl,s,tpsi1,tpsi2,tpsi3,area(l),p_I) +...
comp_quad(Bl,xl(1),yl(1),s,cphi,alpha*R,area(l),p_I,w_I,p_T_12);
\end{Verbatim}

The function \verb+triangle_quad+ estimates $I_{\ell,\ell}^{i,j}$, while 
\verb+comp_quad+ computes numerically the value of $J_{\ell}^{i,j}$. These 
functions use
pre-built data from the file \verb+data.mat+: the first one employs the 
variables \verb+tpsi1+, \verb+tpsi2+ and \verb+tpsi3+, and the second one 
\verb+cphi+, \verb+p_I+, \verb+w_I+ and \verb+p_T_12+. Implementation details 
can be found in appendixes \ref{sec:identical} and \ref{sec:comp}, 
respectively. 
The output of both \verb+triangle_quad+ and \verb+comp_quad+ are 3 by 3 
matrices, such that:
$$ \verb+triangle_quad+_{ij} \approx I^{i,j}_{\ell,\ell}, \quad 
\verb+comp_quad+_{ij} \approx 2J^{i,j}_{\ell} . $$

\subsection{Non-touching elements} \label{ss:non_touching}
The next step is to compute the interactions between $T_\ell$ and all the 
elements $T_m$ whose closure is disjoint $\overline{T_\ell}$ (so that their 
indices are stored in the variable \verb+empty+). In order to do this, we 
calculate and store quadrature points for all the triangles involved in the 
operation as follows:   

\begin{Verbatim}[fontsize=\small]
BBm(:,1:2*ll) = reshape( [ empty_vtx( : , 2:3:3*ll ) -...
  empty_vtx( : , 1:3:3*ll ) , ...
  empty_vtx( : , 3:3:3*ll ) -...
  empty_vtx( : , 2:3:3*ll ) ] , [] , 2)' ;
vl = p_T_6*(Bl') + [ ones(6,1).*xl(1) ones(6,1).*yl(1) ]; 
vm(1:6*ll,:) = reshape( permute( reshape( p_T_6*BBm(:,1:2*ll), ...
 [6 1 2 ll] ) , [1 4 3 2] ) , [ 6*ll 2 ] ) +...
 empty_vtx(: , aux_ind(1:6*ll) )'; 
\end{Verbatim}
The matrix \verb+BBm+ has size $2\times 2*\verb+nt+$, and it contains 
$\verb+nt+$ submatrices of dimension $2\times 2$. The $m$-th submatrix corresponds to the affine transformation that maps $\hat{T}$ into $T_m$. 
The vectors \verb+vl+ and \verb+vm+ contain the coordinates of all quadrature 
points in 
$T_{\ell}$ and $T_m$ for $m \in \verb+empty+$, respectively.

Here, the matrix \verb+BBm+ satisfies 
$$ \verb+BBm(:,2*m-1:2m)'+ \cdot \hat{T} + \texttt{empty\_vtx(:,3*(m-1) + 1)'} 
\mapsto T_{m}, $$

The matrix \verb+p_T_6+ $\in\R^{6 \times 2}$ was provided by the precomputed 
file \verb+data.mat+, and it stores the coordinates of the $6$ quadrature points 
in the reference element $\hat{T}$. In order to compute \verb+vm+, we use three 
nested operations over the $6\times2*\verb+ll+$ matrix 
\verb+p_T_6*BBm(:,1:2*ll)+. To better understand this, suppose we rewrite this 
matrix as follows:
$$\verb+p_T_6*BBm(:,1:2*ll)+ = [ A_1 , A_2 , ... , A_{\verb+ll+}], $$ 
where $A_i$ is a $6\times2$ matrix and $i=1,..,\verb+ll+$. Then, after the 
application of \verb+reshape( permute( reshape( ... '+, we obtain the 
6*\verb+ll+ by 2 matrix $[ A_1 ; A_2 ; ... ; A_{\verb+ll+}]$, which can be used 
as an input in \verb+pdist2+. This trick was taken out from \cite{TrucosMatlab}. 
   
Next, we compute distances from all the quadrature nodes in \verb+vl+ to the 
ones in \verb+vm+, and raise them to the power of $-(2+ 2s)$: 
\begin{Verbatim}[fontsize=\small]
norms(:,1:ll) = reshape(pdist2(vl,vm(1:6*ll,:)),36,[]).^(-2-2*s);
\end{Verbatim}
Thereby, \verb+norms+ is a $36\times\verb+ll+$ matrix such that for $\verb+m+ 
\in \{1,...,\verb+ll+\}$, 
\[
\verb+norms(:,m)+ = 
  \begin{pmatrix} ||\verb+vl(1,:)+ - \texttt{ vm(6*m - 5,:) } ||^{-(2 + 2s)} \\ 
											
\vdots \\
  ||\verb+vl(1,:)+ - \texttt{ vm(6*m, \quad \, :) } ||^{-(2 + 2s)} \\ 
  ||\verb+vl(2,:)+ - \texttt{ vm(6*m - 5,:) } ||^{-(2 + 2s)} \\
\vdots \\ 
 ||\verb+vl(2,:)+ - \texttt{ vm(6*m, \quad \, :) } ||^{-(2 + 2s)}\\
\vdots \\ 
||\verb+vl(6,:)+ - \texttt{ vm(6*m - 5,:) } ||^{-(2 + 2s)}\\
\vdots \\ 
||\verb+vl(6,:)+ - \texttt{ vm(6*m, \quad \, :) } ||^{-(2 + 2s)} \end{pmatrix} ,
\]
where $|| \cdot ||$ denotes the usual euclidean distance in $\mathbb{R}^2$.

At this point, we have collected all the necessary information to compute 
$I_{\ell,m}^{i,j}$ for $\overline{T_{\ell}} \cap \overline{T_m} = \emptyset $ 
and
$i, j$ corresponding to any of the six vertices of these elements. 
We employ the pre-built matrices \verb+phiA+, \verb+phiB+ and \verb+phiD+, that 
contain the values of the nodal basis functions evaluated at the quadrature 
points of $\hat{T}$, multiplied by their respective weights, and stored in an 
appropriate way in order to perform an efficient vectorized operation. Details 
are provided in 
%\hyperref[sec:nontouching]{Appendix \ref*{sec:nontouching}} 
appendixes \ref{sec:nontouching} and \ref{sec:nontouchingdata}. The code 
proceeds: 
\begin{Verbatim}[fontsize=\small]
ML(1:3,1:3,1:ll) =  reshape( phiA*norms(:,1:ll) , 3 , 3 , [] );
ML(1:3,4:6,1:ll) =  reshape( phiB*norms(:,1:ll) , 3 , 3 , [] );
ML(4:6,4:6,1:ll) =  reshape( phiD*norms(:,1:ll) , 3 , 3 , [] );
ML(4:6,1:3,1:ll) =  permute( ML(1:3,4:6,1:ll) , [2 1 3] ) ; 
\end{Verbatim}
So, the matrix \verb+ML+ satisfies 
\[
I_{\ell,m}^{i,j} \approx 4 |T_\ell| |T_m| \, \verb+ML(i,j,m)+ .
\]

The last step to complete the computations for the case $\overline{T_{\ell}} 
\cap \overline{T_m} = \emptyset$ is to add the calculated values in their 
corresponding stiffness matrix entries:
\begin{Verbatim}[fontsize=\small]
for m=1:ll
    order = [nodl t( empty(m) , : )];
    K(order,order) = K(order,order) +...
    ( 8*area(empty(m))*area(l) ).*ML(1:6,1:6,m);
end
\end{Verbatim}
The vector \verb+order+ collects the local indices of the vertices of $T_\ell$ 
and $T_m$, given as explained in Remark \ref{rem:numeracion_ap}. 
Recall that $I_{\ell,m}^{i,j} = I_{m,\ell}^{i,j}$ and that we are summing over 
the elements listed in \verb+empty+. 
In particular, this means that $\ell < m$. We multiply \verb+ML(1:6,1:6,m)+ by 
\verb+8*area(empty(m))*area(l)+ instead of by \verb+4*area(empty(m))*area(l)+ in 
order to avoid carrying the redundant computation of $I_{m,\ell}^{i,j}$.

\subsection{Vertex-touching elements}
In order to compute $I_{\ell,m}^{i,j}$ for the indices $m$ corresponding to 
elements sharing a vertex with $T_\ell$, we use the pre-built variables 
\verb+vpsi1+, \verb+vpsi2+ and \verb+p_cube+ as input in the function 
\verb+vertex_quad+. Let us mention once more that \verb+vpsi1+ and \verb+vpsi2+ 
contain the nodal basis in the reference element $\hat T$ evaluated at 
quadrature points, multiplied by their respective weight and properly stored. 
Moreover, the variable \verb+p_cube+ stores quadrature nodes in the unit cube 
$[0,1]^3$.
Further details about \verb+vertex_quad+ and the auxiliary pre-built data can be 
found in appendixes \ref{sec:vertex} and \ref{sec:vertexdata}, respectively.
We compute the integrals and add the resulting values to \verb+K+ as follows: 
\begin{Verbatim}[fontsize=\small]
for m=vertex  
    nodm = t(m,:);
    nod_com = intersect(nodl, nodm); 
    order = [nod_com nodl(nodl~=nod_com) nodm(nodm~=nod_com)];
    K(order,order) = K(order,order) ... 
    + 2.*vertex_quad(nodl,nodm,nod_com,p,s,vpsi1,vpsi2,...
    area(l),area(m),p_cube); 
end
\end{Verbatim}
Here, we store in \verb+nodm+ the indices of the vertices of $T_m$, whereas 
\verb+nod_com+ dentoes the index of the vertex shared by $T_\ell$ and $T_m$. The 
first entry of \verb+order+ is the index of this common vertex, followed by the 
nodes of $T_\ell$ different from it, and then by the indices of the remaining 
two nodes of $T_m$.
Observe that, unlike the previous case, here there are involved five nodal 
basis, 
%being $\varphi_1$ associated with the common vertex $x_1 = T_{\ell} \cap T_m$, 
so the output of \verb+vertex_quad+ is a 5 by 5 array, such that:
$$
\verb+vertex_quad+_{ij} \approx 
I_{\ell,m}^{i,j}.
$$      

\subsection{Edge-touching elements}
Proceeding similarly, we compute next the case where $\overline{T_{\ell}} \cap 
\overline{T_m}$ is an edge. Now there are only 4 nodal basis functions 
involved, 
and the local numbering is such that the first two nodes correspond to the 
endpoints of the shared edge, the third is the one in $T_\ell$ but not in $T_m$ 
and the last one is the node in $T_m$ but not in $T_\ell$.
Using the pre-built variables \verb+epsi1+, \verb+epsi2+, \verb+epsi3+, 
\verb+epsi4+,\verb+epsi5+ and \verb+p_cube+ as input in \verb+edge_quad+ (see 
appendixes \ref{sec:edge} and \ref{sec:edgedata}), we proceed as in the previous 
case: 
\begin{Verbatim}[fontsize=\small]
for m=edge
    nodm = t(m,:);
    nod_diff = [setdiff(nodl, nodm) setdiff(nodm, nodl)]; 
    order = [ nodl( nodl~=nod_diff(1) ) nod_diff  ];
    K(order,order) = K(order,order) +...
    2.*edge_quad(nodl,nodm,nod_diff,p,s,... 
    epsi1,epsi2,epsi3,epsi4,epsi5,area(l),area(m),p_cube); 
end
\end{Verbatim} 
The indices of the two nodes not shared by $T_\ell$ and $T_m$ are stored in 
\verb+nod_diff+, and \verb+order+ has the nodes ordered as explained in the 
previous paragraph. The output of the function \verb+edge_quad+ is a 4 by 4 
array satisfying
$$ 
\verb+edge_quad+_{ij} \approx 
I_{\ell,m}^{i,j}.
$$  

\subsection{Discrete solution}
Once the main loop is concluded, the stiffness matrix \verb+K+ and the right 
hand side vector \verb+b+ have been computed, and thus it is possible to 
calculate the FE solution \verb+uh+ of the system \eqref{eq:sistema}:
\begin{Verbatim}[fontsize=\small]
uh(nf) = ( K(nf,nf)\b(nf) )./cns; % Solving linear system
\end{Verbatim}
The entries of \verb+K+ and \verb+b+ needed are only the ones corresponding to 
free nodes. The nodes belonging to $\partial \W$ and to the auxiliary domain 
$B\setminus\W$ are excluded, as the discrete solution \verb+uh+ is set to vanish 
on them.

Finally, \verb+uh+ is displayed, and the auxiliary domain is excluded from the 
representation:
\begin{Verbatim}[fontsize=\small]
trimesh(t(1:nt-nt_aux , :), p(1,:),p(2,:),uh);
\end{Verbatim}

\section{Numerical Experiments}
\label{sec:numex}
In order to illustrate the performance of the code, in this section we show the results we obtained in an example problem. 
Explicit solutions for \eqref{eq:fraccionario} are scarce, but it is possible to obtain a family of them if $\W$ is a ball. 
Other numerical experiments carried with this code can be found in \cite{AcostaBorthagaray} and in \cite{BdPM} (for the eigenvalue problem in several domains).

According to the theory  given in \cite{AcostaBorthagaray, BdPM} convergence in the energy norm is expected to occur with order $\frac12$ with respect to the mesh size parameter $h$, or equivalently, of order $-\frac1{2n}$ with respect to the number of degrees of freedom. Moreover, using duality arguments, it is expected to have order of convergence $s+\frac12$ (resp. $- \frac{s+1/2}{n}$) for $0<s \le 1/2$ and $1$ (resp. $-\frac1n$) for $s>1/2$ in the $L^2(\W)$-norm with respect to $h$ (resp. number of degrees of freedom). 

We first construct non-trivial solutions for \eqref{eq:fraccionario} if $\W$ is a ball. Consider the  Jacobi polynomials $P_k^{(\alpha, \beta)} \colon [-1,1] \to \R,$ given by 
{\small{
$$
P_k^{(\alpha, \beta)}(z) =  \frac{\Gamma (\alpha+k+1)}{k!\,\Gamma (\alpha+\beta+k+1)} 
\sum_{m=0}^k {k\choose m} \frac{\Gamma (\alpha + \beta + k + m + 1)}{\Gamma (\alpha + m + 1)} 
\left(\frac{z-1}{2}\right)^m,
$$
}}
 and the weight function  $\w^s:\R^n\to \R,$
 $$\w^s(x) = (1 - \|x\|^2)_+ ^s.$$
 
In \cite[Theorem 3]{DydaKuznetzovKwasnicki} it is shown how to construct explicit  
eigenfunctions for an operator closely
related to the FL by using $P_k^{(s,n/2-1)}$. To be more precise, the authors prove the following result.
\begin{theorem} \label{teo:example} Let $B(0,1)\subset \R^n$ the unitary ball. For $s \in (0,1)$ and $k \in \N$, define 
$$\lambda_{k,s} = \frac{2^{2s} \, \Gamma(1+s+k) \Gamma\left(\frac n2+s+k\right)}{k! \, \Gamma\left(\frac n2+k\right)} $$
and $p_k^{(s)} \colon \rn \to \R$,
$$
p_k^{(s)} (x) = P_k^{(s, \, n/2-1)} ( 2 \|x\|^2 - 1) \chi_{B(0,1)}(x).
$$
Then the following equation holds  
$$
(-\Delta)^s \left(\w^s p_k^{(s)}(x)\right) = \lambda_{k,s} \,  p_k^{(s)}(x) \ \mbox{in } B(0,1).
$$ 
\end{theorem}

A family of explicit solutions is available by using this theorem. 
As a first example, we analyze the solution with $k=0$. This gives a right hand side 
equal to a constant. Namely, consider
\begin{equation} 
\left\lbrace
  \begin{array}{rl}
       (-\Delta)^s u = 1  & \mbox{ in } B(0,1) \subset \R^2, \\
      u = 0 & \mbox{ in }B(0,1)^c . \\
      \end{array}
    \right.
\label{eq:numerical_facil} \end{equation}
We have run the code for a wide range of parameters $s$, while keeping the radius of the auxiliary ball $B$ equal to $1.1$. Orders of convergence in the $L^2$ and energy norm\footnote{A discussion about how to compute errors in the energy norm can be found in \cite{AcostaBorthagaray}.} are shown in Table \ref{tab:ejemplo_1}; these results are in accordance with the theory.
\begin{table}[htbp] 
\caption{Computational rates of convergence for problem \eqref{eq:numerical_facil} with respect to the mesh size, measured in the $L^2(\W)$ and energy norms.}  \label{tab:ejemplo_1}
\footnotesize\centering 
		\begin{tabular}{|c| c| c |} \hline
		Value of $s$ & Order in $L^2(\W)$ & Order in $\widetilde H^s(\W)$ \\ \hline
		$0.1$ & $0.621$  &  $0.500$ \\
		$0.2$ & $0.721$  &  $0.496$ \\
		$0.3$ & $0.804$  &  $0.492$ \\
		$0.4$ & $0.880$  &  $0.491$ \\
		$0.5$ & $0.947$  &  $0.492$ \\
		$0.6$ & $1.003$  &  $0.496$ \\
		$0.7$ & $1.046$  &  $0.501$ \\
		$0.8$ & $1.059$  &  $0.494$ \\
		$0.9$ & $0.999$  &  $0.467$ \\
 \hline	 
\end{tabular}
\end{table}

As a second example we illustrate, in Table \ref{tab:ejemplo_bola},  that in problem \eqref{eq:numerical_facil} the radius $R$ of the auxiliary ball $B$ does not substantially 
affect the error of the scheme. This suggests that it is preferable to maintain the exterior ball's  radius as small as possible.
 Since in this problem the domain $\W$ is itself a ball, for comparison, we also included the output of the code without resorting to the exterior
ball (the row corresponding to $R = 1.0$).
The table clearly shows that the CPU time grows linearly with respect to the number of elements $N_{\tilde{\mathcal{T}}}-N_{\mathcal{T}}$ used in the auxiliary domain.
Taking into account that the final size of the linear system \eqref{eq:sistema} involved in each case is the same,
the computational cost is, essentially, increased only during the assembling routine. Since 
considering an auxiliary domain involves only the computation of the interaction between inner and outer nodes,
a linear behavior of the type described above is clearly expected.

\begin{table}[ht]
\centering
\caption{The $L^2(\Omega)$ and $\widetilde H^s(\W)$ errors for different values of $R$ in problem \eqref{eq:numerical_facil} with $s = 0.5$. In all the cases 
we are using a fixed and regular triangulation $\mathcal{T}$ of $\W$, with $N_{\mathcal{T}} = 4228$. The computations were performed
with {\it MATLAB}\textsuperscript{\textregistered} version 2015a in Windows 10, Intel i7 Processor, RAM 8Gb.}
\label{tab:ejemplo_bola}
\begin{tabular}{|c|c|c|c|c|}
\hline
$R$ & $N_{\tilde{\mathcal{T}}}$ & CPU time (sec.) & Error in $\| \cdot \|_{L^2(\Omega)}$ & Error in $ \| \cdot \|_{\widetilde H^s(\W)}$ \\ \hline
$1.0$ & $4228$ & 80.3 & $0.0164$ & $0.1314$ \\ 
$1.1$ & $4980$ & 100.7 & $0.0167$ & $0.1345$ \\ %\hline
$1.4$ & $8218$ & 206.6  & $0.0167$ & $0.1351$ \\ %\hline
$1.7$ & $12370$ & 344.7 & $0.0167$ & $0.1352$ \\ %\hline
$2.0$ & $17170$ & 511.9 & $0.0167$ & $0.1354$ \\ \hline
\end{tabular}
\end{table}

As a third example we return to the setting of Theorem \ref{teo:example}. We consider 
$k=2$ and compute the order of convergence in $L^2(\W)$ for $s=0.25$ and $s=0.75$. 
We summarize our numerical results in Figure \ref{fig:resultados}. These are in accordance with the predicted rates of 
convergence. Finally, in Figure \ref{fig:grafico} the FE solution, for $s=0.75$ and $k=2$, computed
with a mesh of about 14000 triangles is displayed.

\begin{figure}[ht]
$\begin{array}{cc}
	\includegraphics[width=0.45\textwidth]{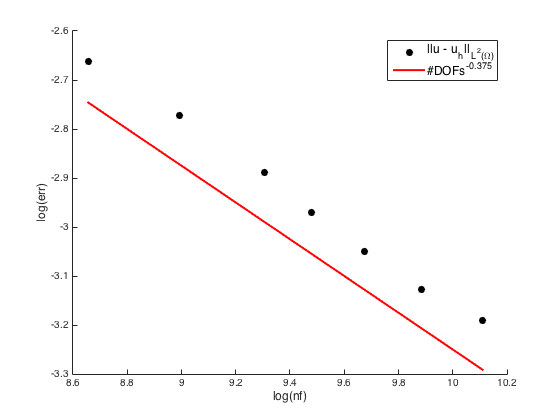} &	\includegraphics[width=0.45\textwidth]{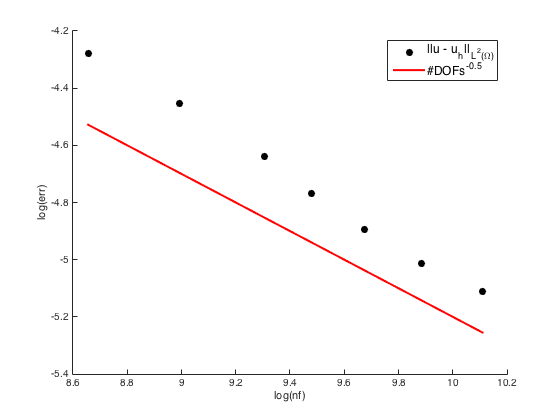} \\
\end{array}$
\caption{Computational rate of convergence in the $L^2(\W)$-norm for the problem with solution given 
by Theorem \ref{teo:example}, for $k = 2$. The left panel corresponds to $s=0.25$ and the right to $s=0.75$. The 
\emph{asymptotic} rate for $s = 0.25$ is $\approx(\#\mbox{DOFs})^{-3/8}$, whereas for $s = 0.75$ it is $\approx(\#\mbox{DOFs})^{-1/2}$, 
in agreement with theory.
} \label{fig:resultados}
\end{figure}

\begin{figure}[ht]
	\includegraphics[width=0.55\textwidth]{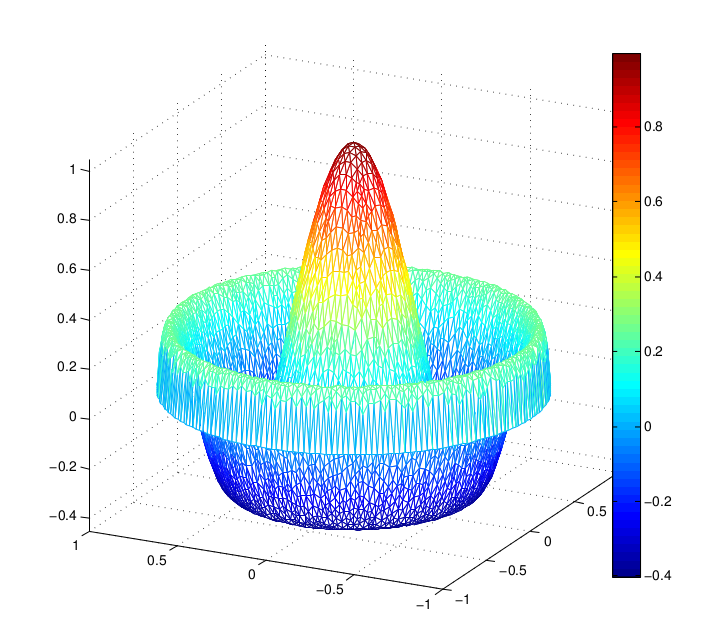}
\caption{FE solution with
a mesh containing about 14000 triangles. With $s=0.75$, we use  $f(x)=\lambda_{2,0.75} \,  p_2^{(0.75)}(x)$
as a source term (see Theorem  \ref{teo:example}).
}
\label{fig:grafico}
\end{figure}
Finally, we would like to mention just a few more facts: our numerical experiments suggest that the condition number of $K$ behaves  like $\sim N_{\mathcal{T}}^{s}$
while over the 99\% of the CPU time is devoted to the assembly routine. Actually, the expected complexity for assembling $K$   
is quadratic in the number of elements, and this seems to be the case in our tests.
 \end{section}

%%%%%%%%%%%%%%%%%%%%%%%%%%%%%%%%%%%%%%%%%%%%%%%%%%%%%%%%%%%%%%%%%%%%%%%%%%%%%%%%
%

\appendix

\section{Quadrature rules} \label{sec:quadrature}

Here we give details about how to compute the integrals $I_{\ell,m}^{i,j}$ and $J_\ell^{i,j}$ (see Section \ref{sec:preliminar}). In order to cope with $I_{\ell,m}^{i,j}$, we proceed according to whether $\overline{T_\ell}\cap \overline{T_m}$ is empty, a vertex, an edge or an element. Recall that $I_{\ell,m} ^{i,j}= I_{m,\ell}^{i,j}$, so that we may assume $\ell \le m$.

Consider two elements $T_\ell$ and $T_m$ such that $\text{supp}(\phii_i),\text{supp}(\phii_j) \cap (T_\ell \cup T_m) \neq \emptyset$.
Observe that if one of this intersections is empty, then $I_{\ell,m}^{i,j} = 0$. Moreover,
 it could be possible that one of the elements is disjoint with the support of both $\phii_i$ and $\phii_j$, provided the other element intersects both supports and  $I_{\ell,m}^{i,j}\neq 0$.

We are going to consider the reference element
\[
\hat{T} = \{ \x = (\hat{x}_1, \hat{x}_2) \colon 0 \le \hat{x}_1 \le 1, \, 0 \le \hat{x}_2 \le \hat{x}_1 \},
\]
whose vertices are 
\[ \x^{(1)} = \left( \begin{array}{c} 0 \\ 0 \end{array} \right), \quad 
\x^{(2)}  = \left( \begin{array}{c} 1 \\ 0 \end{array} \right), \quad 
\x^{(3)} = \left( \begin{array}{c} 1 \\ 1 \end{array} \right).
\]
The basis functions on $\hat{T}$ are, obviously,
\[
\hat{\phii}_1(\x) = 1 - \hat{x}_1,  \quad 
\hat{\phii}_2(\x) = \hat{x}_1- \hat{x}_2, \quad
\hat{\phii}_3(\x) = \hat{x}_2.
\]

\begin{figure}[htbp]
\centering
\begin{tikzpicture}[scale=1.2]
	\draw (-1, -1) -- (0, 0); 
	\draw (0, 0) -- (1.5, -1);
	\draw (1.5, -1) -- (-1, -1);
	\draw (1.5,-1) -- (3,0);
	\draw (3,0) -- (2,1);
	\draw(2,1) -- (1.5,-1); 	
	\node at (0.1, -0.6) {$T_\ell$};
	\node at (2.2, 0) {$T_m$};
	\fill (1.5,-1) circle (2pt);
	\node [below] at (1.5,-1){$1$};	
	\fill (0,0) circle (2pt);
	\node [above] at (0,0){$2$};
	\fill (-1,-1) circle (2pt);
	\node [below] at (-1,-1){$3$};
	\fill (3,0) circle (2pt);
	\node [right] at (3,0){$4$};
	\fill (2,1) circle (2pt);
	\node [above] at (2,1){$5$};
\end{tikzpicture}
\begin{tikzpicture}[scale=1.2]
	\draw (2, 1) -- (0, 0);
	\draw (0, 0) -- (1.5, -1);
	\draw (1.5,-1) -- (3,0);
	\draw (3,0) -- (2,1);
	\draw(2,1) -- (1.5,-1); 	
	\node at (1, 0) {$T_\ell$};
	\node at (2.2, 0) {$T_m$};
	\fill (1.5,-1) circle (2pt);
	\node [below] at (1.5,-1){$1$};	
	\fill (0,0) circle (2pt);
	\node [above] at (0,0){$3$};
	\fill (3,0) circle (2pt);
	\node [right] at (3,0){$4$};
	\fill (2,1) circle (2pt);
	\node [above] at (2,1){$2$};
	\node at (-1,0) {};
\end{tikzpicture}
\caption{Local numbering for elements with a vertex and an edge in common.}
\label{fig:un_punto}
\end{figure}
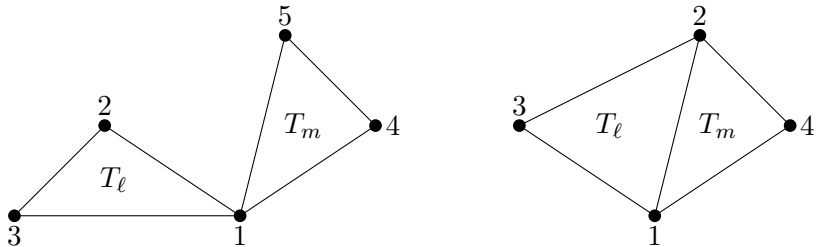

\begin{remark} \label{rem:numeracion_ap}
Given two elements $T_\ell$ and $T_m$, we provide a local numbering in the following way. If $T_\ell$ and $T_m$ are disjoint, we set the first three nodes to be the nodes of $T_\ell$ and the following three nodes to be the ones of $T_m$. Else, we set the first node(s) to be the ones in the intersection, then we insert the remaining node(s) of $T_\ell$ and finally the one(s) of $T_m$ 
(see Figure \ref{fig:un_punto}). For simplicity of notation, when computing $I_{\ell,m}^{i,j}$ and $J_\ell^{i,j}$, we assume that $i,j$ denote the local numbering of the basis functions involved; for example, if $T_\ell$ and $T_m$ share only a vertex, then $1\le i,j\le 5$.
\end{remark}

Consider the affine mappings
\begin{align*}
& \chi_\ell : \hat{T} \to T_\ell, & \chi_\ell(\x) & = B_\ell \x + x_\ell^{(1)} , \\
& \chi_m 		: \hat{T} \to T_m, 		& \chi_m(\x)	&	=  B_m \x + x_m^{(1)} ,
\end{align*}
where the matrices $B_\ell$ and $B_m$ are such that $\x^{(2)}$ (resp. $\x^{(3)}$) is mapped respectively to the second (resp. third) node of $T_\ell$ and $T_m$ in the local numbering defined above. 
Then, it is clear that
{\small{
\begin{equation}  \label{eq:Ilm}
\begin{split}
I^{i,j}_{\ell, m} & = 4 |T_\ell| |T_m| \int_{\hat{T}} \int_{\hat{T}} 
\frac{(\phii_i(\chi_\ell (\x)) - \phii_i(\chi_m (\y)))(\phii_j(\chi_\ell (\x))-\phii_j(\chi_m (\y)))}
{|\chi_\ell (\x) - \chi_m(\y)|^{2+2s}} \, d\x \, d\y = \\
& = 4 |T_\ell| |T_m| 
\iiiint_{\hat{T}\times\hat{T}} F_{ij}(\x_1, \x_2, \y_1, \y_2 ) \, d\x_1 \, d\x_2 \, d\y_1 \, d\y_2.
\end{split}
\end{equation}
}}

We discuss how to compute $I_{\ell,m}^{i,j}$ depending on the relative position of $T_\ell$ and $T_m$, and afterwards we tackle the computation of $J_\ell^{i,j}$.

%%%%%%%%%%%%%%%%%%%%%%%%%%%%%%%%%%%%%%%%%%%%%%%%%%%%%%%%%%%%%%%%%

\subsection{Non-touching elements}
\label{sec:nontouching}

This is the simplest case, since the integrand $F_{ij}$ in \eqref{eq:Ilm} is not singular. Recall that 
$$I^{i,j}_{\ell,m}  = \int_{T_\ell}\int_{T_m} \frac{(\phii_i(x) - \phii_i(y))(\phii_j(x)-\phii_j(y))}{|x-y|^{2+2s}} \, dx dy, \quad 1\le \ell, m \le N_{\mathcal{\tilde T}}. $$
Splitting the numerator in the integrand, we obtain
\begin{align*}
I^{i,j}_{\ell,m}   = & \int_{T_\ell}\int_{T_m} \frac{\phii_i(x) \phii_j(x) }{|x-y|^{2+2s}} \, dx dy  + \int_{T_\ell}\int_{T_m} \frac{\phii_i(y) \phii_j(y) }{|x-y|^{2+2s}}\, dx dy \\ 
- & \int_{T_\ell}\int_{T_m} \frac{\phii_i(x) \phii_j(y) }{|x-y|^{2+2s}} \, dx dy  - \int_{T_\ell}\int_{T_m} \frac{\phii_i(y) \phii_j(x) }{|x-y|^{2+2s}}\, dx dy .
\end{align*}

Note that all the integrands depend on $\ell$ and $m$ only through their denominators.
Since $\phii_i(x)=0$ if $i=1,2,3$ and $x\in T_m$ or if $i=4,5,6$ and $x\in T_{\ell}$, given two indices $i,j$, only one of the four integrals above is not null.
Thus, we may divide the 36 interactions between the 6 basis functions involved into four 3 by 3 blocks, and write the local matrix \verb+ML+ as: 
\begin{equation}
\verb+ML+ = \left( \begin{array}{cc}
A_{\ell,m} & B_{\ell,m} \\ C_{\ell,m} & D_{\ell,m} \end{array} \right),
\label{eq:ML}
\end{equation}
where 
%\peligro{sacar los $\ell, m$ de estos bloques para alivianar la notacion?} \verificar{ Si es cierto que se simplifica la notacion, pero me parece que si los sacamos aca un poco mas abajo los tenemso que volver a poner cuando uso el conjunto $\mathcal{S}$. Ahora estoy dudando de que es lo que convendria}
%
{\small{
\begin{align*}
& A^{i,j}_{\ell,m} = \int_{T_\ell}\int_{T_m} \frac{\phii_i(x) \phii_j(x) }{|x-y|^{2+2s}}\, dx dy, & B^{i,j}_{\ell,m} = -\int_{T_\ell}\int_{T_m} \frac{\phii_i(x) \phii_{j+3}(y) }{|x-y|^{2+2s}} \, dx dy\\
& C^{i,j}_{\ell,m} =-\int_{T_\ell}\int_{T_m} \frac{\phii_{i+3}(y) \phii_{j}(x) }{|x-y|^{2+2s}} \, dx dy, & D^{i,j}_{\ell,m} = \int_{T_\ell}\int_{T_m} \frac{\phii_{i+3}(y) \phii_{j+3}(y) }{|x-y|^{2+2s}}  \, dx dy .
\end{align*}
}}
We use two nested Gaussian quadrature rules to estimate these integrals. These have $6$ quadrature nodes each, making a total of $36$ quadrature points. 
Let us denote by $p_k$ and $w_k$ ($k=1,\ldots,6$) the quadrature nodes and weights in $\hat{T}$, respectively.
%So, for instance, if we want to compute $A_{\ell,m}$, let be $x_1,...,x_6$ the quadrature points over $\hat{T}$ and $w_1,...,w_6$ their respective weights such that, if $f: \hat{T} \rightarrow  \mathbb{R}$ is a regular function, then
%$$ 
%\int_{\hat{T}}fdx \approx \sum_{k=1}^6 w_{k}f(x_k) .
%$$ 
Changing variables we obtain 
$$
A^{i,j}_{\ell,m} =4|T_{\ell}||T_{m}| \int_{\hat{T}}\int_{\hat{T}} \frac{\hat{\phii}_i(x) \hat{\phii}_j(x) }{|\chi_{\ell}(x)-\chi_{m}(y)|^{2+2s}}\, dx dy,
$$
and applying the quadrature rule twice, we derive: 
\begin{equation}
A^{i,j}_{\ell,m} \approx 4 |T_\ell| |T_m|\sum_{q=1}^6 \sum_{k=1}^6 \frac{w_q \, w_k \, \hat{\phii}_i(p_k) \hat{\phii}_j(p_k) }{|\chi_{\ell}(p_k)-\chi_{m}(p_q)|^{2+2s}} .
\label{eq:quad_a}
\end{equation}
Note that the right hand side summands only depend on  $i$ and $j$ through their numerators, and on $\ell$ and $m$ through their denominators. As our goal is to compute the whole block $A_{\ell,m}$ as efficiently as possible, 
we set the following definitions:
\begin{itemize}
\item  The matrix $\Phi^A \in \R^9 \times \R^{36}$ stores the numerators involved in \eqref{eq:quad_a}, corresponding to the $9$ pairs of basis functions and the $36$ pairs of quadrature nodes, respectively. Namely, 
\begin{equation} \label{eq:phiA}
\Phi^A_{ij} = \hat{\varphi}_{[i-1]_{3}+1}(p_{\lceil \frac{j}{6}\rceil}) \hat{\varphi}_{\lceil \frac{i}{3} \rceil }(p_{\lceil \frac{j}{6}\rceil})w_{[j-1]_{6} + 1}w_{\lceil \frac{j}{6}\rceil}, 
\end{equation}
where $[m]_k$ denotes $m$ modulo $k$ and $\lceil \cdot \rceil$ is the ceiling function. Let us make this definition more explicit. 
The matrix $\Phi^{A}$ may be divided in $6$ blocks, 
$$\Phi^{A}=(\Phi^{A_1} \ldots \Phi^{A_6}),$$
where $\Phi^{A_k}$ is a $6 \times 9$ matrix:    
{\small{
$$\Phi^{A_k}=\begin{pmatrix}
\hat{\varphi}_1(p_k)\hat{\varphi}_1(p_k)w_{k} w_{1}& \hat{\varphi}_1(p_k)\hat{\varphi}_1(p_k)w_{k}w_{2} & \ldots & \hat{\varphi}_1(p_k)\hat{\varphi}_1(p_k)w_{k}w_{6} \\
\hat{\varphi}_2(p_k)\hat{\varphi}_1(p_k)w_{k} w_{1}& \hat{\varphi}_2(p_k)\hat{\varphi}_1(p_k)w_{k}w_{2} & \ldots & \hat{\varphi}_2(p_k)\hat{\varphi}_1(p_k)w_{k}w_{6}  \\
\hat{\varphi}_3(p_k)\hat{\varphi}_1(p_k)w_{k} w_{1}& \hat{\varphi}_3(p_k)\hat{\varphi}_1(p_k)w_{k}w_{2} & \ldots & \hat{\varphi}_3(p_k)\hat{\varphi}_1(p_k)w_{k}w_{6}\\[0.3em]
\hat{\varphi}_1(p_k)\hat{\varphi}_2(p_k)w_{k} w_{1}& \hat{\varphi}_1(p_k)\hat{\varphi}_2(p_k)w_{k}w_{2} & \ldots & \hat{\varphi}_1(p_k)\hat{\varphi}_2(p_k)w_{k}w_{6} \\
\hat{\varphi}_2(p_k)\hat{\varphi}_2(p_k)w_{k} w_{1}& \hat{\varphi}_2(p_k)\hat{\varphi}_2(p_k)w_{k}w_{2} & \ldots & \hat{\varphi}_2(p_k)\hat{\varphi}_2(p_k)w_{k}w_{6} \\
\hat{\varphi}_3(p_k)\hat{\varphi}_2(p_k)w_{k} w_{1}& \hat{\varphi}_3(p_k)\hat{\varphi}_2(p_k)w_{k}w_{2} & \ldots & \hat{\varphi}_3(p_k)\hat{\varphi}_2(p_k)w_{k}w_{6}\\[0.3em]
\hat{\varphi}_1(p_k)\hat{\varphi}_3(p_k)w_{k} w_{1}& \hat{\varphi}_1(p_k)\hat{\varphi}_3(p_k)w_{k}w_{2} & \ldots & \hat{\varphi}_1(p_k)\hat{\varphi}_3(p_k)w_{k}w_{6} \\
\hat{\varphi}_2(p_k)\hat{\varphi}_3(p_k)w_{k} w_{1}& \hat{\varphi}_2(p_k)\hat{\varphi}_3(p_k)w_{k}w_{2} & \ldots & \hat{\varphi}_2(p_k)\hat{\varphi}_3(p_k)w_{k}w_{6} \\
\hat{\varphi}_3(p_k)\hat{\varphi}_3(p_k)w_{k} w_{1}& \hat{\varphi}_3(p_k)\hat{\varphi}_3(p_k)w_{k}w_{2} & \ldots & \hat{\varphi}_3(p_k)\hat{\varphi}_3(p_k)w_{k}w_{6}\\
\end{pmatrix} .
$$
}}

\item 
The variable $d^m \in \R^{36}$ is a vector storing the distances between all the quadrature nodes involved:
\begin{equation} \label{eq:dm}
d^m_{k} = \left|\chi_{\ell}(p_{[k-1]_{6}+1})-\chi_{m}(p_{\lceil \frac{k}{6} \rceil}) \right|^{-(2+2s)} . 
\end{equation}
Namely, the vector $d^{m}$ can be written as: 
 $$d^{m}=\begin{pmatrix} \left|\chi_{\ell}(p_{1})-\chi_{m}(p_{1}) \right|^{-(2+2s)} \\ 
											\vdots \\
  \left|\chi_{\ell}(p_{6})-\chi_{m}(p_{1}) \right|^{-(2+2s)} \\ 
  \left|\chi_{\ell}(p_{1})-\chi_{m}(p_{2}) \right|^{-(2+2s)} \\
\vdots \\ 
 \left|\chi_{\ell}(p_{6})-\chi_{m}(p_{2}) \right|^{-(2+2s)}\\
\vdots \\ 
\vdots \\ 
\left|\chi_{\ell}(p_{1})-\chi_{m}(p_{6}) \right|^{-(2+2s)}\\
\vdots \\ 
\left|\chi_{\ell}(p_{6})-\chi_{m}(p_{6}) \right|^{-(2+2s)} \end{pmatrix} .$$\\
\end{itemize}

With these two variables in hand, the computation of the integrals $A^{i j}$ may be done in a vectorized mode. Defining $\hat{A}_{\ell , m} := \Phi^{A} \cdot d^m$, we obtain: 
\begin{align*}
A^{[i-1]_3+1, \lceil \frac{i}{3} \rceil}_{\ell , m} & \approx 4|T_{\ell}||T_{m}|\hat{A}^{i}_{\ell , m} \\ 
& = 4 |T_\ell| |T_m|\sum_{q} \sum_{k} w_q w_k \frac{\hat{\phii}_{[i-1]_3+1}(p_k) \hat{\phii}_{\lceil \frac{i}{3} \rceil}(p_k) }{|\chi_{\ell}(p_k)-\chi_{m}(p_q)|^{2+2s}}, \ i \in \{1,...,9\} .
\end{align*}
Equivalently, using {\it MATLAB}\textsuperscript{\textregistered} notation: 
$$ 
A_{\ell , m} \approx  4|T_{\ell}||T_{m}| \, \verb+reshape(+\hat{A}_{\ell , m} \verb+, 3 , 3)+. 
$$
%\peligro{No entend\'i esta frase: Given that $\Phi^{A}$ do not change along the execution, we only need to compute it once. } \verificar{Quise poner que es una cosa constante, y no depende de nada del problema solamente de la cuadrtura, igual abajo ya lo explicaste mucho mejor} \peligro{Hay que cambiar $n$ por 36 en las expresiones de mas abajo}   
 
We apply the same ideas to computate the remaining blocks in \eqref{eq:ML}. We define: 
\begin{itemize}
\item  a $9\times 36$ matrix $\Phi^B$, such that
 $$
\Phi^B_{ij} = \hat{\varphi}_{[i-1]_{3}+1}(p_{\lceil \frac{j}{n}\rceil}) \hat{\varphi}_{\lceil \frac{i}{3} \rceil +3}(p_{[j-1]_n+1})w_{[j-1]_{n} + 1}w_{\lceil \frac{j}{n}\rceil} ,
$$
\item a $9\times 36$ matrix $\Phi^D$, such that
 $$
\Phi^D_{ij} = \hat{\varphi}_{[i-1]_{3}+4}(p_{[j-1]_n+1}) \hat{\varphi}_{\lceil \frac{i}{3} \rceil + 3}(p_{[j-1]_n+1})w_{[j-1]_{n} + 1}w_{\lceil \frac{j}{n}\rceil} .
$$
\end{itemize}   
Then, considering
\begin{align*}
\hat{B}_{\ell , m} & := \Phi^{B} \cdot d^m, \\
\hat{D}_{\ell , m} & := \Phi^{D} \cdot d^m,
\end{align*}
we just need to multiply
\begin{eqnarray*}
& B_{\ell , m} \approx  4|T_{\ell}||T_{m}|\verb+reshape(+\hat{B}_{\ell , m} \verb+, 3 , 3)+, \\
& D_{\ell , m} \approx  4|T_{\ell}||T_{m}|\verb+reshape(+\hat{D}_{\ell , m} \verb+, 3 , 3)+ . 
\end{eqnarray*}

Is simple to verify that $C_{\ell,m} = B'_{\ell,m}$, so that there is no need to make additional operations to compute the block $C_{\ell, m}$. 
Moreover, let us emphasize that the matrices $\Phi^{A}$, $\Phi^{B}$ and $\Phi^{D}$ depend on the quadrature rule employed, but not on the elements under consideration; these are precomputed and stored in \verb+data.mat+. We refer to Section \ref{sec:nontouchingdata} for details on how this is done. 
%So, given that we can pre-compute $\Phi^{A}$, $\Phi^{B}$ and $\Phi^{D}$, fixed $\ell$ and $m$, we just need to calculate $d^m$. So we can estimate $I_{\ell,m}$ as follow:
However, in the main loop, the vector $d^m$ needs to be calculated for every $1\le \ell \le m \le {N}_{\hat{\mathcal T}}$. 

We obtain a matrix \verb+ML+ as follows:
\small{
$$
\verb+ML+ \approx 4|T_{\ell}||T_{m}| \left( \begin{array}{cc}
\verb+reshape(+\Phi^{A} \cdot d^m \verb+, 3 , 3)+ & \verb+reshape(+\Phi^{B} \cdot d^m \verb+, 3 , 3)+ \\ \verb+reshape(+\Phi^{B} \cdot d^m \verb+, 3 , 3)'+ & \verb+reshape(+\Phi^{D} \cdot d^m \verb+, 3 , 3)+ \end{array} \right).
$$
}
%Nevertheless, there is an efficient way to compute the local matrices \verb+ML+ for a given $\ell$ and several values of $m \in \{1,...,N_{\mathcal{\tilde T}}\}$.
%\peligro{No entiendo esto: se usa en el c\'odigo? Con qu\'e conjunto $\mathcal S$?} \verificar{ $\mathcal{S}$ viene a ser el conjunto de indices que en el codigo principal estan en el vector empty, ahi estan los triangulos que tiene interseccion vacia con $\ell$ y que todavia no fueron recorridos. Como la cuenta esta sirve para calcular rapido la cuadratura contra cualquier subconjunto de elementos lo escribi generico porque me parecio mas claro. Ahora que lo veo si es muy confuso.. por el momento agregue una aclaracion en parentesis. }
%
In addition, this vectorized approach gives us an efficient way to compute $I_{\ell,m}$ for several values of $m \in \{1,...,N_{\mathcal{\tilde T}}\}$ at once. Indeed, suppose that want to calculate $I_{\ell,m}$ for $m \in \mathcal{S} \subseteq \{1,...,N_{\mathcal{\tilde T}}\}$ (along the execution of the main code, $\mathcal{S}$ would contain the indices listed in \verb+empty+). It is possible to compute $\hat{A}_{\ell,m}$, $\hat{B}_{\ell,m}$ and $\hat{D}_{\ell,m}$ for all $m \in \mathcal{S}$ using vectorized operations as follows: 
$$ \left( \hat{A}_{\ell,m_1},..., \hat{A}_{\ell,m_{\#\mathcal{S} } } \right) = \Phi^{A}\cdot \left( d^{m_1},...,d^{m_{\#\mathcal{S} }} \right), $$
$$ \left( \hat{B}_{\ell,m_1},..., \hat{B}_{\ell,m_{\#\mathcal{S} } } \right) = \Phi^{B}\cdot \left( d^{m_1},...,d^{m_{\#\mathcal{S} }} \right), $$
$$ \left( \hat{D}_{\ell,m_1},..., \hat{D}_{\ell,m_{\#\mathcal{S} } } \right) = \Phi^{D}\cdot \left( d^{m_1},...,d^{m_{\#\mathcal{S} }} \right). $$

Observe that, fixed $\ell$ and $\mathcal{S}$, the distances between interpolation points of the involved triangles are all the necessary information to obtain the estimation of the matrix \verb+ML+ (given by \eqref{eq:ML}), for $m \in \mathcal{S}$. 

In order to perform an efficient computation of $\left( d^{m_1},...,d^{m_{\#\mathcal{S} }} \right)$, we use the Matlab function \verb+pdist2+ in the following way: 
$$
 \left( d^{m_1},...,d^{m_{\#\mathcal{S} }} \right)= \verb+reshape( pdist2(+ X^{\ell} ,  \left( \begin{array}{c} X^{m_1} \\  X^{m_2} \\ \vdots \\ X^{m_{\#\mathcal{S}}} \end{array} \right)  \verb+), + n^2  \verb+, [] )+ ^{(-1-s)} .
 $$
Here, the vectors $X^m$ are given by 
$$ X^{m} := \left( \begin{array}{c} \chi_m(p_1) \\  \vdots \\ \chi_m(p_6) \end{array} \right) .$$ 
The computation of the matrix \verb+ML+ is carried in the main code, and it is implemented in {Subsection \ref{ss:non_touching}.

%%%%%%%%%%%%%%%%%%%%%%%%%%%%%%%%%%%%%%%%%%%%%%%%%%%
\subsection{Vertex-touching elements}
\label{sec:vertex}
%%%%%%%%%%%%%%%%%%%%%%%%%%%%%%%%%%%%%%%%%%%%%%%%%%%
In case $\overline{T_\ell} \cap \overline{T_m}$ consists of a vertex, define $\z = (\x, \y)$, identify $\z$ with a vector in $\mathbb{R}^4$, 
 and split the domain of integration in \eqref{eq:Ilm} 
into two components $D_1$ and $D_2,$ where
\begin{align*}
D_1 = \{ \z 
%=(\hat{z}_1,\hat{z}_2,\hat{z}_3,\hat{z}_4) 
\colon 0 \le \hat{z}_1 \le 1, \,	0 \le \hat{z}_2 \le \hat{z}_1, \,
0 \le \hat{z}_3 \le \hat{z}_1, \, 0 \le \hat{z}_4 \le \hat{z}_3 \}, \\
D_2 = \{ \z 
%=(\hat{z}_1,\hat{z}_2,\hat{z}_3,\hat{z}_4) 
\colon 0 \le \hat{z}_3 \le 1, \,	0 \le \hat{z}_4 \le \hat{z}_3, \,
0 \le \hat{z}_1 \le \hat{z}_3, \, 0 \le \hat{z}_2 \le \hat{z}_1 \}.
\end{align*}

Let $\xi \in [0,1]$ and $\eta = (\eta_1, \eta_2, \eta_3) \in [0,1]^3.$ We consider the mappings $T_h \colon [0,1]\times[0,1]^3 \to D_h, \ h =1,2$,
\begin{align*}
T_1(\xi, \eta) = 
\left( \begin{array}{c}
\xi \\ \xi\eta_1 \\ \xi\eta_2 \\ \xi\eta_2\eta_3 \end{array} \right), 
\quad
T_2(\xi, \eta) = 
\left( \begin{array}{c}
\xi\eta_2 \\ \xi\eta_2\eta_3 \\ \xi \\ \xi\eta_1
\end{array} \right),
\end{align*}
having Jacobian determinants $|JT_1| = \xi^3 \eta_2=|JT_2|.$

We perform the calculations in detail only on $D_1$. Observe that if $i=1$, which corresponds to the vertex in common between $T_\ell$ and $T_m$, then
\[
\phii_i(\chi_\ell (\xi, \xi\eta_1)) - \phii_i(\chi_m (\xi\eta_2, \xi\eta_2\eta_3)) = - \xi (1 - \eta_2) .
\]
Meanwhile, if the subindex $i$ equals $2$ or $3$, it corresponds to one of the other two vertices of $T_\ell$. Therefore, in those cases $\phii_i(\chi_m (\xi\eta_2, \xi\eta_2\eta_3))=0$, and
\begin{align*}
\phii_2(\chi_\ell (\xi, \xi\eta_1)) = &  \xi (1 - \eta_1), \\
\phii_3(\chi_\ell (\xi, \xi\eta_1)) = & \xi \eta_1 .
\end{align*}
Analogously, if $i \in\{4, 5\}$, then $\phii_i(\chi_\ell (\xi, \xi\eta_1))=0$ and so
\begin{align*}
\phii_4(\chi_m (\xi\eta_2, \xi\eta_2\eta_3)) = &  -\xi \eta_2 (1 - \eta_3), \\
\phii_5(\chi_m (\xi\eta_2, \xi\eta_2\eta_3)) = &  -\xi \eta_2 \eta_3 .
\end{align*}

Thus, defining the functions $\psi^{(1)}_k  \colon [0,1]^3 \to \R$ ($k \in \{1,\ldots, 5\}$),
\begin{align*}
\psi^{(1)}_1 (\eta) & = \eta_2 - 1, & \psi^{(1)}_2 (\eta) & = 1- \eta_1, 
& \psi^{(1)}_3 (\eta) & = \eta_1, \\
\psi^{(1)}_4 (\eta) & = -\eta_2  (1 - \eta_3), & \psi^{(1)}_5 (\eta) & = -\eta_2\eta_3, &
\end{align*}
we may write
\begin{align*}
\int_{D_1} F_{ij}(\z) \, d\z
= & \int_{[0,1]} \int_{[0,1]^3} \frac{\psi^{(1)}_i(\eta) \psi^{(1)}_j(\eta)}
{\left|
B_\ell \left( \begin{array}{c} \xi \\ \xi\eta_1 \end{array} \right) - 
B_m \left( \begin{array}{c}\xi\eta_2 \\ \xi\eta_2\eta_3 \end{array} \right)
\right|^{2+2s}} \, \xi^5 \eta_2 \, d\eta \, d\xi \\
= & \left(\int_0^1 \xi^{3-2s} d\xi \right)
\left(\int_{[0,1]^3} \frac{\psi^{(1)}_i(\eta) \psi^{(1)}_j(\eta)}
{\left| d^{(1)}(\eta) \right|^{2+2s}} \,
\eta_2 \, d\eta \right) \\
= & \frac{1}{4-2s} \left(\int_{[0,1]^3} \frac{\psi^{(1)}_i(\eta) \psi^{(1)}_j(\eta)}
{\left| d^{(1)}(\eta) \right|^{2+2s}} \,
\eta_2 \, d\eta \right),
\end{align*}
where we have defined the function
\[ 
d^{(1)}(\eta) = B_\ell \left( \begin{array}{c} 1 \\ \eta_1 \end{array} \right) - 
B_m \left( \begin{array}{c}\eta_2 \\ \eta_2\eta_3 \end{array} \right).
\]
Observe that  
in the first line of last equation (or equivalently, in \eqref{eq:Ilm}), the integrand is singular at the origin. The key point in the identity above is that the singularity of the integral is explicitly computed. The function $d^{(1)}$ is not zero on $[0,1]^3$, and
therefore the last integral involves a regular integrand that is easily estimated by means of a Gaussian quadrature rule.

In a similar fashion, the integrals over $D_2$ take the form
\begin{align*}
\int_{D_2} F_{ij}(\z) \, d\z
= & \frac{1}{4-2s} \left(\int_{[0,1]^3} \frac{\psi^{(2)}_i(\eta) \psi^{(2)}_j(\eta)}
{\left| d^{(2)}(\eta) \right|^{2+2s}} \,
\eta_2 \, d\eta \right) ,
\end{align*}
where
\begin{align*}
\psi^{(2)}_1(\eta) & = 1 - \eta_2, & \psi^{(2)}_2(\eta) & = \eta_2 (1 - \eta_3), &
\psi^{(2)}_3(\eta) & = \eta_2\eta_3, \\
\psi^{(2)}_4(\eta) & = \eta_1 - 1, & \psi^{(2)}_5(\eta) & = - \eta_1, &
\end{align*}
and
\[
d^{(2)}(\eta)= B_\ell \left( \begin{array}{c} \eta_2 \\ \eta_2\eta_3 \end{array} \right) 
- B_m  \left( \begin{array}{c} 1 \\ \eta_1 \end{array} \right)  .
\]

Based on the previous analysis, we describe the function \verb+vertex_quad+.  Let $p_1,...,p_n \in [0,1]^3$ be a set of quadrature points and $w_1,...,w_n$ their respective weights. In the code we present, we work with three nested three-point quadrature rules on $[0,1]$, making a total of $27$ quadrature nodes in the unit cube. The data necessary to use this quadrature is supplied in the file \verb+data.mat+, and in Appendix \ref{sec:points}.
%, such that if $f: [0,1]^3 \rightarrow  \mathbb{R}$ is a regular function, then
%$$  \int_{[0,1]^3}fdx \approx \sum_k w_{k}f(p_k) $$
%We consider the case over $Dh$
%\peligro{Unificar indice dominios auxiliares: se usa $D_h$, $D_k$...}

Set $h \in \{ 1,2 \} $. Then, applying the mentioned quadrature rule in the cube, 
$$
\int_{[0,1]^3} \frac{\psi^{(h)}_i(\eta) \psi^{(h)}_j(\eta)}
{\left| d^{(h)}(\eta) \right|^{2+2s}} \,
\eta_2 \, d\eta   \approx  \sum_{k=1}^{27}  w_k \frac{\psi^{(h)}_i(p_k) \psi^{(h)}_j(p_k)}
{\left| d^{(h)}(p_k) \right|^{2+2s}} \,
p_{k,2} ,
$$ 
where $p_{k,2}$ denotes the second coordinate of the point $p_k$.
The right hand side only depends on $\ell$ and $m$ through $d^{(h)}$. So, in order to compute $I_{\ell,m}^{i,j}$ using vectorized operations, we define the following variables, in analogy to \eqref{eq:phiA} and \eqref{eq:dm}:
\begin{itemize}
\item 
A $25\times 27$ matrix $\Psi^h$ satisfying
$$
\Psi^h_{ij} = w_j  \, \psi^{(h)}_{[i-1]_5 + 1}(p_j) \psi^{(h)}_{\lceil \frac{i}{5} \rceil}(p_j) \, p_{j,2}.$$   
\item 
A vector $d^h \in \R^{27}$, such that $$d^h_{k} = \left| d^{(h)}(p_k) \right|^{2+2s} .$$
\end{itemize}   

Then, defining $\hat{I}_{\ell , m} := \Psi^1 \cdot d^1 +  \Psi^2 \cdot d^2$, we obtain 
\begin{align*}
 I^{[i-1]_5+1, \lceil \frac{i}{5} \rceil}_{\ell , m} & \approx \frac{4|T_{\ell}||T_{m}|}{4 - 2s}\hat{I}^{i}_{\ell , m}  \\
&=  \sum_{h=1}^2 \sum_{k=1}^{27}  w_k \frac{\psi^{(h)}_{[i-1]_5+1}(p_k) \psi^{(h)}_{\lceil \frac{i}{5} \rceil}(p_k)}
{\left| d^{(h)}(p_k) \right|^{2+2s}}, \quad i \in \{1,...,25\}.
\end{align*}
Equivalently, using {\it MATLAB}\textsuperscript{\textregistered} notation: 
$$ I_{\ell , m} \approx  \frac{4|T_{\ell}||T_{m}|}{4 - 2s}\verb+reshape(+\hat{I}_{\ell , m} \verb+, 5 , 5)+ .$$
Given that the matrices $\Psi^{1}$ and $\Psi^{2}$ do not change along the execution, we only need to compute them once. These are precomputed and provided on the data file; explicit information regarding its entries is available on Appendix \ref{sec:vertexdata}.
 
So, the function \verb+vertex_quad+ computes the previous quadrature rule in the following way:
\begin{Verbatim}[fontsize=\small]
function ML = vertex_quad (nodl,nodm,sh_nod,p,s,psi1,psi2,areal,aream,p_c)
xm = p(1, nodm);
ym = p(2, nodm);
xl = p(1, nodl);
yl = p(2, nodl);
x = p_c(:,1);
y = p_c(:,2);
z = p_c(:,3);
local_l = find(nodl==sh_nod);
nsh_l = find(nodl~=sh_nod);
nsh_m = find(nodm~=sh_nod);
p_c = [xl(local_l), yl(local_l)];
Bl = [xl(nsh_l(1))-p_c(1) xl(nsh_l(2))-xl(nsh_l(1));
      yl(nsh_l(1))-p_c(2) yl(nsh_l(2))-yl(nsh_l(1))];
Bm = [xm(nsh_m(1))-p_c(1) xm(nsh_m(2))-xm(nsh_m(1));
      ym(nsh_m(1))-p_c(2) ym(nsh_m(2))-ym(nsh_m(1))];  
ML = ( 4*areal*aream/(4-2*s) ).*reshape(...
     psi1*( sum( ([ones(length(x),1)  x]*(Bl')...
    	 - [y , y.*z]*(Bm') ).^2, 2 ).^(-1-s) ) +...
     psi2*( sum( ([ones(length(x),1)  x]*(Bm')...
    	 - [y , y.*z]*(Bl') ).^2, 2 ).^(-1-s) )...
     , 5 , 5);
end
\end{Verbatim}
In the code above, \verb+nodl+ and \verb+nodm+ are the vertex indices of $T_\ell$ and $T_m$ respectively, \verb+sh_nod+ is the index of the shared node, \verb+p+ is an array that contains all the vertex coordinates, \verb+areal+ and \verb+aream+ denote $|T_\ell|$ and $|T_m|$ respectively, \verb+s+ is $s$, and \verb+p_c+ contains the coordinates of the quadrature points on $[0,1]^3$.  This last variable is gathered form \verb+data.mat+, where it is stored as \verb+p_cube+ (see Appendix \ref*{sec:points}). 
In addition, \verb+Bl+ and \verb+Bm+ play the role of $B_{\ell}$ and $B_{m}$, and \verb+psi1+ and \verb+psi2+ are $\Psi^{1}$ and $\Psi^{2}$ respectively. 
As we mentioned, 
\verb+psi1+ and \verb+psi2+ have been pre-computed and stored on \verb+data.mat+ as \verb+vpsi1+ and \verb+vpsi2+ respectively (see Appendix \ref*{sec:vertexdata}).  

The output of \verb+vertex_quad+ is a $6\times 6$ matrix \verb+ML+  that satisfies $\verb+ML(i,j)+ \approx  I_{\ell, m}^{i,j}$.    

%%%%%%%%%%%%%%%%%%%%%%%%%%%%%%%%%%%%%%%%%%%%%%%%%%%
\subsection{Edge-touching elements} 
\label{sec:edge}
%%%%%%%%%%%%%%%%%%%%%%%%%%%%%%%%%%%%%%%%%%%%%%%%%%%
In this case, the parametrization of the elements we are considering is such that both $\chi_\ell$ and $\chi_m$ map $[0,1]\times\{0\}$ to the common edge between $T_\ell$ and $T_m$. Therefore, if we consider $\z = (\y_1 - \x_1, \y_2, \x_2)$, the singularity of the integrand is localized at $\z = 0$:
\[
I^{i,j}_{\ell,m} = 4 |T_\ell||T_m| \int_0^1 \int_{-\x_1}^{1-\x_1} \int_{0}^{\z_1+\x_1} \int_{0}^{\x_1} 
F_{ij}(\x_1, \z_3, \x_1+\z_1, \z_2) \, d\z \, d\x_1 .
\]

We decompose the domain of integration as $\cup_{k=1}^{5} D_k$, where 
\begin{align*}
D_1  = \{ (\x_1,\z) \colon & -1 \le \hat{z}_1 \le 0, \, 0 \le \hat{z}_2 \le 1+\hat{z}_1, \\
&  0 \le \z_3 \le \z_2-\z_1, \, \z_2-\z_1 \le \x_1 \le 1 \}, \\
D_2  = \{ (\x_1,\z) \colon & -1 \le \hat{z}_1 \le 0, \,	0 \le \hat{z}_2 \le 1+\hat{z}_1, \\
&  \z_2-\z_1\le \z_3\le 1, \, \z_3 \le \x_1 \le 1 \}, \\
D_3  = \{ (\x_1,\z) \colon & 0 \le \hat{z}_1 \le 1, \,	0 \le \hat{z}_2 \le \hat{z}_1, \\
& 0 \le \z_3\le 1-\z_1, \, \z_3 \le \x_1 \le 1-\z_1 \}, \\
D_4 = \{ (\x_1,\z) \colon & 0 \le \hat{z}_1 \le 1, \,	\z_1 \le \hat{z}_2 \le 1, \\
&  0 \le \z_3\le \z_2-\z_1, \, \z_2-\z_1 \le \x_1 \le 1 - \z_1 \}, \\
D_5  = \{ (\x_1,\z) \colon & 0 \le \hat{z}_1 \le 1, \,	\z_1 \le \hat{z}_2 \le 1, \\
&  \z_2-\z_1 \le \z_3\le 1-\z_1, \, \z_3 \le \x_1 \le 1 - \z_1 \} .
\end{align*}
Consider the mappings $T_k \colon [0,1]\times[0,1]^3 \to D_k$ ($k \in \{1,\ldots, 5\}$), 
\begin{align*} 
& T_1 \left( \begin{array}{c} \xi \\ \eta \end{array} \right)
= \left( \begin{array}{c}  \xi \\ -\xi\eta_1\eta_2 \\ \xi\eta_1(1-\eta_2) \\ \xi\eta_1\eta_3 \end{array} \right), 
& T_2\left( \begin{array}{c} \xi \\ \eta \end{array} \right) 
= \left( \begin{array}{c} \xi \\ -\xi\eta_1\eta_2\eta_3 \\ \xi\eta_1\eta_2(1-\eta_3) \\ \xi\eta_1 \end{array} \right), \\
& T_3\left( \begin{array}{c} \xi \\ \eta \end{array} \right)
= \left( \begin{array}{c}  \xi(1-\eta_1\eta_2) \\ \xi\eta_1\eta_2 \\ \xi\eta_1\eta_2\eta_3 \\ \xi\eta_1(1-\eta_2) \end{array} \right),  
& T_4\left( \begin{array}{c} \xi \\ \eta \end{array} \right)
= \left( \begin{array}{c}  \xi(1-\eta_1\eta_2\eta_3) \\  \xi\eta_1\eta_2\eta_3 \\ \xi\eta_1 \\ \xi\eta_1\eta_2(1-\eta_3) \end{array} \right), \\
& T_5\left( \begin{array}{c} \xi \\ \eta \end{array} \right)
= \left( \begin{array}{c}  \xi(1-\eta_1\eta_2\eta_3) \\ \xi\eta_1\eta_2\eta_3 \\ \xi\eta_1\eta_2 \\ \xi\eta_1(1-\eta_2\eta_3) \end{array} \right), &   
\end{align*}
with Jacobian determinants given by
\[
|JT_1| = \xi^3\eta_1^2, \quad |JT_h| =  \xi^3\eta_1^2\eta_2, \ h\in \{2,\ldots , 5 \}.
\]
Then, over $D_h$ it holds that
\begin{align*}
\int_{D_h} F_{ij} = \frac{1}{4-2s} \int_{[0,1]^3} \frac{\psi^{(h)}_i(\eta) \psi^{(h)}_j(\eta)}
{|d^{(h)}(\eta)|^{2+2s}} \,
J^{(h)}(\eta) \, d\eta ,
\end{align*}
where 
\begin{align*}
\psi_1^{(1)}(\eta) & = - \eta_1\eta_2, & \psi_2^{(1)}(\eta) & = \eta_1 (1 - \eta_3), \\
\psi_3^{(1)}(\eta) & = \eta_1\eta_3, & \psi_4^{(1)}(\eta) & = -\eta_1 (1-\eta_2), \\
\psi_1^{(2)}(\eta) & = - \eta_1\eta_2\eta_3, & \psi_2^{(2)}(\eta) & = -\eta_1 (1 - \eta_2), \\ \psi_3^{(2)}(\eta) & = \eta_1, & \psi_4^{(2)}(\eta) & = -\eta_1 \eta_2(1-\eta_3), \\
\psi_1^{(3)}(\eta) & = \eta_1\eta_2, & \psi_2^{(3)}(\eta) & = -\eta_1 (1 - \eta_2\eta_3), \\ \psi_3^{(3)}(\eta) & = \eta_1(1-\eta_2), & \psi_4^{(3)}(\eta) & = -\eta_1 \eta_2\eta_3,\\ 
\psi_1^{(4)}(\eta) & = \eta_1\eta_2\eta_3, & \psi_2^{(4)}(\eta) & = \eta_1 (1 - \eta_2), \\ \psi_3^{(4)}(\eta) & = \eta_1\eta_2(1-\eta_3), & \psi_4^{(4)}(\eta) & = -\eta_1 , \\
\psi_1^{(5)}(\eta) & = \eta_1\eta_2\eta_3, & \psi_2^{(5)}(\eta) & = -\eta_1 (1 - \eta_2), \\ \psi_3^{(5)}(\eta) & = \eta_1(1-\eta_2\eta_3), & \psi_4^{(5)}(\eta) & = -\eta_1 \eta_2.
\end{align*}
Moreover, the functions $d^{(h)}$ are given by
\begin{align*}
d^{(1)}(\eta) & = B_\ell \left( \begin{array}{c} 1 \\ \eta_1\eta_3 \end{array} \right) 
- B_m \left( \begin{array}{c} 1-\eta_1\eta_2 \\ \eta_1(1-\eta_2) \end{array} \right),\\
d^{(2)}(\eta) & = B_\ell \left( \begin{array}{c} 1 \\ \eta_1 \end{array} \right) 
- B_m\left( \begin{array}{c} 1-\eta_1\eta_2\eta_3 \\ \eta_1\eta_2(1-\eta_3) \end{array} \right) ,\\
d^{(3)}(\eta) & = B_\ell \left( \begin{array}{c} 1-\eta_1\eta_2 \\ \eta_1(1-\eta_2) \end{array} \right)
- B_m \left( \begin{array}{c} 1 \\ \eta_1\eta_2\eta_3 \end{array} \right),\\
d^{(4)}(\eta) & = B_\ell \left( \begin{array}{c} 1-\eta_1\eta_2\eta_3 \\ \eta_1\eta_2(1-\eta_3) \end{array} \right)
- B_m \left( \begin{array}{c} 1 \\ \eta_1 \end{array} \right),\\
d^{(5)}(\eta) & = B_\ell \left( \begin{array}{c} 1-\eta_1\eta_2\eta_3 \\ \eta_1(1-\eta_2\eta_3) \end{array} \right)
- B_m \left( \begin{array}{c} 1 \\ \eta_1\eta_2 \end{array} \right),
\end{align*}
and the Jacobians are
\[
J^{(1)}(\eta) =\eta_1^2, \quad J^{(h)}(\eta) =\eta_1^2\eta_2, \ h\in \{2,\ldots , 5 \}.
\]

As in the case of vertex-touching elements, the problem is reduced to computing integrals on the unit cube. Let $p_1,...,p_{27} \in [0,1]^3$ the quadrature points, and $w_1,...,w_{27}$ their respective weights. For $h = 1, \ldots ,5$ we have
$$
\int_{[0,1]^3} \frac{\psi^{(h)}_i(\eta) \psi^{(h)}_j(\eta)}
{\left| d^{(h)}(\eta) \right|^{2+2s}} \,
J^{(h)}(\eta) \, d\eta   \approx  \sum_k  w_k \frac{\psi^{(h)}_i(p_k) \psi^{(h)}_j(p_k)}
{\left| d^{(h)}(p_k) \right|^{2+2s}} \,
J^{(h)}(p_k) .$$ 
Once more, the right hand side only depends on $\ell$ and $m$ through $d^{(h)}$. So, with the purpose of computing $I_{\ell,m}$ efficiently,  we define: 
\begin{itemize}
\item 
A matrix $\Psi^h \in \R^{16\times 27}$, given by 
$$\Psi^h_{ij} = w_j \, \psi_{[i-1]_4+1}(p_j) \psi_{\lceil \frac{i}{4} \rceil}(p_j)J^{(h)}(p_j) .$$   
\item 
A vector $d^h \in \R^{27}$, such that 
$$d^h_{k} = \left| d^{(h)}(p_k) \right|^{2+2s} .$$
\end{itemize}   
Therefore,  considering $\hat{I}_{\ell , m} = \Psi^1 \cdot d^1 + \dots +  \Psi^5 \cdot d^5$, we reach the following relation: 
\begin{align*}
I^{[i-1]_4+1, \lceil \frac{i}{4} \rceil}_{\ell , m} & \approx \frac{4|T_{\ell}||T_{m}|}{4 - 2s} \hat{I}^{i}_{\ell , m}  \\
& = \sum_{h}\sum_k  w_k \frac{\psi^{(h)}_{[i-1]_4+1}(p_k) \psi^{(h)}_{\lceil \frac{i}{4} \rceil}(p_k)}
{\left| d^{(h)}(p_k) \right|^{2+2s}}, \: i \in \{1,...,16\} .
\end{align*}
Using {\it MATLAB}\textsuperscript{\textregistered} notation,
$$ I_{\ell , m} \approx \frac{4|T_{\ell}||T_{m}|}{4 - 2s} \verb+reshape(+\hat{I}_{\ell , m} \verb+, 4 , 4)+ .$$
As before, the matrices $\Psi^{1}$, \ldots, $\Psi^{5}$ do not depend on the elements under consideration, so they are precomputed and provided in \verb+data.mat+, where they are stored as $\verb+epsi1+, \ldots, \verb+epsi5+$, respectively. Details about their calculation are given in Appendix \ref{sec:edgedata}.      

The function \verb+edge_quad+ performs the calculations we have explained in this section.
\begin{Verbatim}[fontsize=\small]
function ML = edge_quad(nodl,nodm,nod_diff,p,s,psi1,psi2,psi3,...
                        psi4,psi5,areal,aream,p_c)
xm = p(1, nodm);
ym = p(2, nodm);
xl = p(1, nodl);
yl = p(2, nodl);
x = p_c(:,1);
y = p_c(:,2);
z = p_c(:,3);
local_l = find(nodl~=nod_diff(1));
nsh_l = find(nodl==nod_diff(1));
nsh_m = find(nodm==nod_diff(2));
P1 = [xl(local_l(1)), yl(local_l(1))];
P2 = [xl(local_l(2)), yl(local_l(2))];
Bl = [P2(1)-P1(1) -P2(1)+xl(nsh_l);
      P2(2)-P1(2) -P2(2)+yl(nsh_l)];
Bm = [P2(1)-P1(1) -P2(1)+xm(nsh_m);
      P2(2)-P1(2) -P2(2)+ym(nsh_m)];
ML = ( 4*areal*aream/(4-2*s) ).*reshape(... 
      psi1*( sum( ([ones(length(x),1)  x.*z]*(Bl')...
	 - [1-x.*y  x.*(1-y)]*(Bm') ).^2, 2 ).^(-1-s) ) +...
      psi2*( sum( ([ones(length(x),1)  x]*(Bl')...
	 - [1-x.*y.*z  x.*y.*(1-z)]*(Bm') ).^2, 2 ).^(-1-s) ) +...
      psi3*( sum( ([(1-x.*y)  x.*(1-y)]*(Bl')...
	 - [ones(length(x),1) x.*y.*z]*(Bm') ).^2, 2 ).^(-1-s) ) +...
      psi4*( sum( ([1-x.*y.*z  x.*y.*(1-z)]*(Bl')...
	 - [ones(length(x),1) x]*(Bm') ).^2, 2 ).^(-1-s) ) +... 
      psi5*( sum( ([1-x.*y.*z  x.*(1-y.*z)]*(Bl')...
	 - [ones(length(x),1) x.*y]*(Bm') ).^2, 2 ).^(-1-s) )...
      , 4 , 4);
end
\end{Verbatim}
Here, \verb+nodl+ and \verb+nodm+ are the indices of the vertices of $T_\ell$ and $T_m$ respectively, \verb+nod_diff+ contains the not-shared-vertex index, \verb+p+ is an array that contains all the vertex coordinates, \verb+areal+ and \verb+aream+ are $|T_\ell|$ and $|T_m|$ respectively, \verb+s+ is $s$, \verb+p_c+ contains the coordinates of the quadrature points on $[0,1]^3$ (stored in \verb+data.mat+, see Appendix \ref{sec:points}), \verb+Bl+ and \verb+Bm+ are $B_{\ell}$ and $B_{m}$, and \verb+psi1+, ..., \verb+psi5+ are $\Psi^{1},\dots,\Psi^{5}$ respectively.
%As before, \verb+psi1+,..., \verb+psi5+ have been pre-computed and stored on \verb+data.mat+ as \verb+epsi1+,..., \verb+epsi5+ respectively (see \hyperref[sec:edgedata]{Appendix \ref*{sec:edgedata}}).  

The output of this function is a $4\times 4$ matrix $\verb+ML+ \approx  I_{\ell, m}$.

%%%%%%%%%%%%%%%%%%%%%%%%%%%%%%%%%%%%%%%%%%%%%%%%%%%
\subsection{Identical elements} 
\label{sec:identical}
%%%%%%%%%%%%%%%%%%%%%%%%%%%%%%%%%%%%%%%%%%%%%%%%%%%
In the same spirit as before, let us consider $\z = \y - \x$, so that
\[
I_{\ell,\ell} = 4 |T_\ell|^2 
\int_0^1 \int_0^{\x_1} \int_{-\x_1}^{1-\x_1} \int_{-\x_2}^{\z_1+\x_1-\x_2}
%\frac{(\hat{\phii}_i(\x)-\hat{\phii}_i(\x+\z))(\hat{\phii}_j(\x)-\hat{\phii}_j(\x+\z))}
%{|B_\ell (\z)|^{2+2s}} 
F_{ij}(\x_1, \x_2, \x_1+\z_1, \x_2+\z_2)
\, d\z_2 \, d\z_1 \, d\x_2 \, d\x_1 .
\]

Let us decompose the integration region into 
\begin{equation} \label{eq:descomposicion}\begin{split}
D_1  = \{ (\x, \z) \colon & -1 \le \z_1 \le 0, \ -1\le \z_2\le\z_1, \\ & -\z_2 \le \x_1 \le 1, \ -\z_2 \le \x_2 \le \x_1 \},\\
D_2 = \{ (\x, \z) \colon & 0 \le \z_1 \le 1, \ \z_1\le \z_2\le 1, \\ & \z_2-\z_1 \le \x_1 \le 1-\z_1, \ 0 \le \x_2 \le \z_1 - \z_2 + \x_1 \},\\
D_3  = \{ (\x, \z) \colon  &-1 \le \z_1 \le 0, \ \z_1\le \z_2\le 0, \\ & -\z_1 \le \x_1 \le 1, \ -\z_2 \le \x_2 \le \x_1+\z_1-\z_2 \},\\
D_4 = \{ (\x, \z) \colon & 0 \le \z_1 \le 1, \ 0\le \z_2\le \z_1, \\ &  0 \le \x_1 \le 1-\z_1, \ 0 \le \x_2 \le \x_1 \},\\
D_5  = \{ (\x, \z) \colon & -1 \le \z_1 \le 0, \ 0\le \z_2\le 1+\z_1, \\ &  \z_2-\z_1 \le \x_1 \le 1, \ 0 \le \x_2 \le \x_1+\z_1-\z_2 \},\\
D_6  = \{ (\x, \z) \colon & 0 \le \z_1 \le 1, \ -1+\z_1 \le \z_2\le 0, \\ & -\z_2 \le \x_1 \le 1-\z_1, \ -\z_2 \le \x_2 \le \x_1 \}.
\end{split}\end{equation}

We begin by considering the first two sets. Making the change of variables $(\x',\z') = (\x, -\z)$ on $D_1$ and $(\x', \z')=(\x+\z,\z)$ on $D_2$, both regions are transformed into
\[
D'_1 = \{(\x',\z') \colon 0\le \z'_1\le 1, \ \z'_1\le\z'_2\le 1, \ \z'_2\le \x'_1 \le 1, \ \z'_2\le \x'_2\le \x'_1 \},
\]
so that
\begin{align*}
4|T_\ell|^2 \int_{D_1\cup D_2} F_{ij}(\x,\x+\z) & =
4|T_\ell|^2 \int_{D'_1} F_{ij}(\x',\x'-\z') +  F_{ij} (\x'-\z', \x') \, d\x' \, d\z' \\
& = 8 |T_\ell|^2 \int_{D'_1}\ F_{ij}(\x',\x'-\z') \, d\x' \, d\z',
\end{align*}
because
\[
F_{ij} (\x', \x'-\z') = 
\frac{ (\hat{\phii}_i (\x') - \hat{\phii}_i (\x'-\z'))(\hat{\phii}_j (\x') - \hat{\phii}_j (\x'-\z'))}
{|B_\ell(\z') |^{2+2s}}
 = F_{ij} (\x'-\z',\x').
\] 
Next, consider the four-dimensional simplex
\[
D = \{ w \colon 0\le w_1 \le 1, \ 0\le w_2 \le w_1, \ 0\le w_3 \le w_2, \ 0\le w_4 \le w_3 \},
\]
the map $T_1 \colon D \to D'_1$,
\[
 \left( \begin{array}{c} \x' \\ \z' \end{array} \right) = 
T_1 \left( \begin{array}{c} w_1 \\ w_2 \\ w_3 \\ w_4  \end{array} \right)  = 
\left( \begin{array}{c} w_1, \\ w_1-w_2+w_3, \\ w_4, \\ w_3 \end{array} \right), \quad |JT_1| = 1, 
\]
and the Duffy-type transform $T: [0,1]^4 \to D$,
\begin{equation} \label{eq:duffyD}
w = T \left( \begin{array}{c} \xi \\ \eta \end{array} \right) = 
\left( \begin{array}{c} \xi, \\ \xi\eta_1, \\ \xi\eta_1\eta_2, \\ \xi\eta_1\eta_2\eta_3 \end{array} \right), \quad |JT| = \xi^3\eta_1^2\eta_2.
\end{equation}

The composition of these two changes of variables allows to write the variables in $F_{ij}$ in terms of $(\xi,\eta)$ in the following way:
\begin{align*}
\x' = \left( \begin{array}{c} \xi \\ \xi (1- \eta_1 + \eta_1 \eta_2) \end{array} \right), \
\z' = \left( \begin{array}{c} \xi\eta_1\eta_2\eta_3 \\ \xi \eta_1 \eta_2 \end{array} \right), \
\x-\z' = \left( \begin{array}{c} \xi (1- \eta_1 \eta_2 \eta_3) \\ \xi (1- \eta_1) \end{array} \right) .
\end{align*}

Observe that
\[
\Lambda_k^{(1)}(\xi,\eta) := \hat{\phii}_k (\x') - \hat{\phii}_k (\x'-\z') = 
\begin{cases}
-\xi \eta_1 \eta_2 \eta_3 & \text{ if } k = 1, \\
-\xi \eta_1 \eta_2 (1-\eta_3) &  \text{ if } k = 2, \\
\xi \eta_1 \eta_2 & \text{ if } k = 3.
\end{cases}
\]
Thus,
\begin{align*}
 4|T_\ell|^2 \int_{D_1\cup D_2} F_{ij}(\x,\x+\z) & = 8|T_\ell|^2 \int_D F_{ij}(w_1,w_1-w_2+w_3,w_4, w_3) \, dw = \\
&= 8|T_\ell|^2 \int_{[0,1]^4} \frac{\Lambda_i^{(1)}(\xi,\eta) \, \Lambda_j^{(1)}(\xi,\eta)}
{\left| B_\ell \left( \begin{array}{c} \xi \eta_1 \eta_2 \eta_3 \\ \xi \eta_1 \eta_2 \end{array} \right) \right|^{2+2s}}
\, \xi^3 \eta_1^2 \eta_2 \, d\xi \, d\eta .
\end{align*}
Finally, as the functions $\Lambda_k^{(1)}$ may be rewritten as $\Lambda_k^{(1)}(\xi,\eta) = \xi \eta_1 \eta_2 \psi_k^{(1)}(\eta_3)$,
where
\[
\psi_1^{(1)}(\eta_3) = -\eta_3, \quad \psi_2^{(1)}(\eta_3) = -(1-\eta_3), \quad \psi_3^{(1)}(\eta_3) = 1 ,
\]
we obtain
\begin{align*}
& 4|T_\ell|^2 \int_{D_1\cup D_2} F_{ij}(\x,\x+\z) = \\
& = 8 |T_\ell|^2 \int_0^1 \xi^{3-2s} d\xi \ \int_0^1 \eta_1^{2-2s} d\eta_1 \
\int_0^1 \eta_2^{1-2s} d\eta_2 \ \int_0^1 \frac{\psi_i^{(1)}(\eta_3)\psi_j^{(1)}(\eta_3)}
{\left|B_\ell \left( \begin{array}{c} \eta_3 \\ 1 \end{array} \right) \right|^{2+2s}} d\eta_3.
\end{align*}

Obviously, the first three integrals above are straightforwardly calculated by hand, and the last one involves a regular integrand, so that it is easily estimated by means of a Gaussian quadrature rule.

It still remains to perform similar calculations on the rest of the sets in \eqref{eq:descomposicion}. Consider the new variables
$(\x',\z') =(\x,-\z)$ on $D_3$, $(\x',\z') =(\x+\z,\z)$ on $D_4$, $(\x',\z') =(\x+\z,\z)$ on $D_5$ and $(\x',\z') =(\x,-\z)$ on $D_6$, 
so that
\begin{align*}
4|T_\ell|^2 \int_{D_3\cup D_4} F_{ij}(\x,\x+\z) & = 8 |T_\ell|^2 \int_{D'_2}\ F_{ij}(\x',\x'-\z') \, d\x' \, d\z', \\
4|T_\ell|^2 \int_{D_5\cup D_6} F_{ij}(\x,\x+\z) & = 8 |T_\ell|^2 \int_{D'_3}\ F_{ij}(\x',\x'-\z') \, d\x' \, d\z',
\end{align*}
where
\begin{align*}
D'_2 & = \{(\x',\z') \colon 0\le \z'_1\le 1, \ 0\le\z'_2\le \z'_1, \ \z'_1\le \x'_1 \le 1, \ \z'_2\le \x'_2\le \x'_1 -\z'_1 + \z'_2 \}, \\
D'_3 & = \{(\x',\z') \colon -1\le \z'_1\le 0, \ 0\le\z'_2\le 1+\z'_1, \ \z'_2\le \x'_1 \le 1+\z'_1, \ \z'_2\le \x'_2\le \x'_1 \}.
\end{align*}
These domains are transformed into $[0,1]^4$ by the respective composition of the transformations $T_h \colon D \to D'_h$ ($h=1,2$)
\[
T_2 \left( \begin{array}{c} w_1 \\ w_2 \\ w_3 \\ w_4 \end{array} \right) =  \left( \begin{array}{c} w_1 \\ w_2-w_3+w_4 \\ w_3 \\ w_4 \end{array} \right), \quad
T_3 \left( \begin{array}{c} w_1 \\ w_2 \\ w_3 \\ w_4 \end{array} \right) =  \left( \begin{array}{c} w_1-w_4 \\ w_2-w_4 \\ -w_4 \\ w_3- w_4 \end{array} \right),
\]
and the Duffy transformation \eqref{eq:duffyD}. Simple calculations lead finally to
{\small{
\begin{align*}
4|T_\ell|^2 \int_{D_3\cup D_4} \! \! F_{ij}(\x,\x+\z) & = \frac{8 |T_\ell|^2}{(4-2s)(3-2s)(2-2s)} \int_0^1 \frac{\psi_i^{(2)}(\eta_3)\psi_j^{(2)}(\eta_3)}
{\left|B_\ell \left( \begin{array}{c} 1 \\ \eta_3 \end{array} \right) \right|^{2+2s}} d\eta_3, \\
4|T_\ell|^2 \int_{D_5\cup D_6} \! \! F_{ij}(\x,\x+\z) & = \frac{8 |T_\ell|^2}{(4-2s)(3-2s)(2-2s)} \int_0^1 \frac{\psi_i^{(3)}(\eta_3)\psi_j^{(3)}(\eta_3)}
{\left|B_\ell \left( \begin{array}{c} \eta_3 \\ 1-\eta_3 \end{array} \right) \right|^{2+2s}} d\eta_3,
\end{align*}
}}
where
\begin{align*}
\psi_1^{(2)} (\eta_3) & = - 1, & \psi_2^{(2)} (\eta_3) & = 1- \eta_3,  & \psi_3^{(2)} (\eta_3) & = \eta_3, \\
\psi_1^{(3)} (\eta_3) & = \eta_3, & \psi_2^{(3)} (\eta_3) & = - 1,  & \psi_3^{(3)} (\eta_3) & = 1-\eta_3.
\end{align*}

For the sake of simplicity of notation, we write
\begin{align*}
d^{(1)} (x) & := \left|B_\ell \left( \begin{array}{c} x \\ 1 \end{array} \right) \right|^{2+2s},  & d^{(2)} (x) & := \left|B_\ell \left( \begin{array}{c} 1 \\ x \end{array} \right) \right|^{2+2s},\\ 
 d^{(3)} (x) & := \left|B_\ell \left( \begin{array}{c} x \\ 1 - x \end{array} \right) \right|^{2+2s} . &
\end{align*}

In order to estimate the integrals in the unit interval, we use a $9$ point Gaussian quadrature rule. Let $p_1, \ldots ,p_9 \in [0,1]$ the quadrature points, and $w_1,...,w_9$ their respective weights. Considering the integrals over the domains $D'_h$ ($h  \in  \{1, 2, 3 \}$), we may write
$$
\int_0^1 \frac{\psi_i^{(h)}(\eta)\psi_j^{(h)}(\eta)}
{d^{(h)}(\eta)} d\eta
  \approx  \sum_{k=1}^9  w_k \frac{\psi_i^{(h)}(p_k)\psi_j^{(h)}(p_k)}
{d^{(h)}(p_k)} .
$$ 

As before, we take advantage of the fact that the integrand only depends on $\ell$ through its denominator. We define:
\begin{itemize}
\item 
A $9\times9$ matrix $\Psi^h$, such that 
$$\Psi^h_{ij} = w_j \, \psi_{[i-1]_3 + 1}(p_j) \psi_{\lceil \frac{i}{3} \rceil}(p_j) \, J^{(h)}(p_j) .$$   
\item 
A vector $d^h \in \R^9$, given by 
$$d^h_{k} = d^{(h)}(p_k) .$$
\end{itemize}   
Setting 
$\hat{I}_{\ell , m} := \Psi^1 \cdot d^1 + \Psi^2 \cdot d^2 +  \Psi^3 \cdot d^3$, we obtain, for  $ i \in \{1,...,9\}$,
\begin{align*} 
I^{[i-1]_3 + 1, \lceil \frac{i}{3} \rceil}_{\ell , m} & \approx \frac{8 |T_\ell|^2}{(4-2s)(3-2s)(2-2s)}\hat{I}^{i}_{\ell , m}  \\ 
& = \frac{8 |T_\ell|^2}{(4-2s)(3-2s)(2-2s)}\sum_{h=1}^3\sum_{k=1}^9  w_k \frac{\psi_{[i-1]_3 + 1}^{(h)}(p_k)\psi_{\lceil \frac{i}{3} \rceil}^{(h)}(p_k)}
{d^{(h)}(p_k)} ,
\end{align*}
or in {\it MATLAB}\textsuperscript{\textregistered} notation, 
$$ 
I_{\ell , m} \approx \frac{8 |T_\ell|^2}{(4-2s)(3-2s)(2-2s)} \, \verb+reshape(+\hat{I}_{\ell , m} \verb+, 3 , 3)+ .
$$
The matrices $\Psi^{1}$, $\Psi^{2}$ and $\Psi^{3}$ are supplied by \verb+data.mat+, where they are respectively saved as \verb+tpsi1+, \verb+tpsi2+ and \verb+tpsi3+.

The code of the function \verb+triangle_quad+ is as follows.
\begin{Verbatim}[fontsize=\small]
function ML = triangle_quad(Bl,s,psi1,psi2,psi3,areal,p_I)
ML = ( 8*areal*areal/((4-2*s)*(3-2*s)*(2-2*s)) ).*reshape(...
  psi1*( (  sum( (Bl*[p_I'; ones(1,length(p_I))]).^2 ).^(-1-s) )' ) +...
  psi2*( (  sum( (Bl*[ones(1,length(p_I)) ; p_I']).^2 ).^(-1-s) )' ) + ...
  psi3*( (  sum( (Bl*[p_I' ; p_I' - ones(1,length(p_I))]).^2 ).^(-1-s) )' ) ...
  , 3 , 3);
end
\end{Verbatim}
The matrix \verb+Bl+ plays the role of $B_{\ell}$, \verb+s+ is $s$, \verb+areal+ is $|T_{ \ell}|$, and \verb+p_I+ contains the values of the quadrature points in $[0,1]$. The latter are stored in \verb+data.mat+ under the same name, see Appendix \ref{sec:points}. The matrices $\Psi^{1}$, $\Psi^{2}$ and $\Psi^{3}$ are respectively saved as \verb+psi1+, \verb+psi2+ and \verb+psi3+.
%Again, \verb+psi1+, \verb+psi2+ and \verb+psi3+ have been pre-computed and stored in \verb+data.mat+ as \verb+tpsi1+, \verb+tpsi2+ and \verb+tpsi3+ respectively (see \hyperref[sec:identicaldata]{appendix \ref*{sec:identicaldata} }).   

The output \verb+ML+ of this function is a $3\times3$ matrix, such that: $\verb+ML+ \approx  I_{\ell, \ell}$.

%%%%%%%%%%%%%%%%%%%%%%%%%%%%%%%%%%%%%%%%%%%%%%%%%%%
\subsection{Complement}
\label{sec:comp}
%%%%%%%%%%%%%%%%%%%%%%%%%%%%%%%%%%%%%%%%%%%%%%%%%%%
Recall that we are assuming that the domain $\W$ is contained in a ball $B=B(0,R)$.  Here we are considering the interaction of two basis functions $\phii_i$, $\phii_j$ such that $\supp (\phii_i) \cap \supp (\phii_j) = T_\ell$, over the region $T_\ell \times B^c$. Namely, we aim to compute
\begin{align*} \label{eq:complemento}
J_{\ell} & = \int_{T_\ell} \int_{B^c} \frac{\phii_i(x) \phii_j(x)}{|x-y|^{2+2s}} \, dy dx =
\int_{T_\ell} \phii_i(x) \phii_j(x) \psi(x) \, dx  \\
& = 2 |T_\ell| \int_{\hat{T}} \hat{\phii_i}(\x) \hat{\phii_j}(\x) \psi(\chi_\ell(\x)) \, d\x,
\end{align*}
where 
\[
\psi(x) = \int_{B^c} \frac{1}{|x-y|^{2+2s}} \, dy .
\]

The integral above may be calculated by a Gauss quadrature rule in the reference element $\hat T$, provided 
that the values of $\psi$ at the quadrature points are computed.

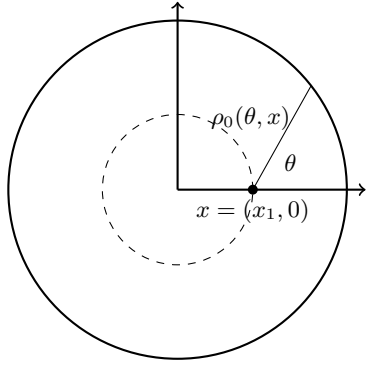
\begin{figure}[ht]
\centering
\begin{tikzpicture}[scale=1.25]
	\draw [thick] (0, 0) circle (1.8cm); 	\draw [dashed] (0, 0) circle (0.8cm); 
	\draw [thick, ->] (0, 0) -- (2, 0); 	\draw [thick, ->] (0, 0) -- (0, 2);
	\draw (0.8, 0) -- (1.42,1.1061);	\fill (0.8,0) circle (1.5pt);
	\node [below] at (0.8,0){\footnotesize{$x = (x_1,0)$}}; 
	\node [left] at (1.3,0.76){\footnotesize{$\rho_0(\theta,x)$}};
	\node [above] at (1.2,0.1){\footnotesize{$\theta$}};
\end{tikzpicture}
		\caption{Computing $\psi(x)$ in a point of $B=B(0,R)$. 
		Due to the symmetry, the value of $\psi$ is the
		same along the dashed circle, hence we may assume that  $x=(x_1,0)$ and $0\le x_1 <R$.
		For any $0\le \theta \le \pi$ , the function $\rho_0$ is given by  
		$\rho_0(\theta,x) = -x_1 \cos \theta + \sqrt{R^2 - x_1^2 \sin^2 \theta}$.}
	\label{fig:afuera}
\end{figure}

Observe that the function $\psi$ is radial (see Figure \ref{fig:afuera}) and therefore it suffices to estimate it on points of the form $x=(x_1, 0)$, 
where $x_1>0$. For a fixed point $x$ and given $\theta \in [0, 2\pi]$, let $\rho_0(\theta)$ be the distance between $x$ 
and the intersection of the ray starting from $x$ with angle $\theta$ with respect to the horizontal axis. Then, it is simple to verify that
\[
\rho_0(\theta,x) = -x_1 \cos \theta + \sqrt{R^2 - x_1^2 \sin^2 \theta} ,
\]
and therefore, integrating in polar coordinates,
\[
\psi(x) = \frac{1}{2s} \int_0^{2\pi} \frac{1}{\rho_0(\theta,x)^{2s}} \, d\theta .
\]

In order to compute $J_{\ell}$ we perform two nested quadrature rules: one over $\hat{T}$ and, for each quadrature point $p_k$ in $\hat{T}$, another one to estimate $\psi(p_k)$ over $[0,2\pi]$. We apply a $12$ point quadrature formula over $\hat{T}$ and a $9$ point one on $[0,2 \pi ]$.
Let $p_1,\ldots,p_{12} \in \hat{T}$, $\theta_1,\ldots,\theta_9 \in [0,2\pi]$ be these quadrature nodes, and $w_1,\ldots,w_{12}$, $W_1,\ldots,W_9$ their respective weights. Applying the rules we obtain
$$J_{\ell}
  \approx  \frac{|T_{\ell}|}{s}\sum_{k=1}^{12}  w_k \hat{\varphi}_i(p_k)\hat{\varphi}_j(p_k) \sum_{q=1}^9  \frac{W_q}{\rho_0(\theta_q,\chi_\ell(p_k))^{2s}}  .$$

In the same fashion as for the other computations, we write the previous expression  as the product of a pre-computed matrix (that only depends on the choice of the quadrature rules) times a vector that depends on the elements under consideration. Indeed, we define:
\begin{itemize}
\item 
A matrix $\Phi \in \R^{9\times12}$, such that 
$$\Phi_{ij} = w_j \, \hat{\varphi}_{[i-1]_3 + 1}(p_j)  \hat{\varphi}_{\lceil \frac{i}{3} \rceil}(p_j).$$
\item 
A vector $\rho \in \R^{12}$,  such that
$$\rho_{k} = \sum_q  \frac{W_q}{\rho_0(\theta_q,\chi_\ell(p_k))^{2s}} .$$
\end{itemize}   
Upon  defining $\hat{J}_{\ell } := \Phi \cdot \rho $, we obtain
$$ 
J^{[i-1]_3 + 1, \lceil \frac{i}{3} \rceil}_{\ell } \approx \frac{|T_{\ell}|}{s}\hat{J}^{i}_{\ell }, \: i \in \{1,...,9\} . $$
Using {\it MATLAB}\textsuperscript{\textregistered} notation, the above identity may be written as
$$ J_{\ell } \approx \frac{|T_{\ell}|}{s} \, \verb+reshape(+\hat{J}_{\ell } \verb+, 3 , 3)+ .$$
%Given that $\Phi$ does not change along the execution, we only need to compute it once.     

The function \verb+comp_quad+ perform the previous computations.
\begin{Verbatim}[fontsize=\small]
function ML = comp_quad(Bl, x0, y0, s , phi , R, areal , p_I , w_I , p_T)
x = (Bl*p_T')' + [x0.*ones(length(p_T),1) , y0.*ones(length(p_T),1)];
aux = x(:,1)*cos(2*pi*p_I') + x(:,2)*sin(2*pi*p_I');
weight = ( ( -aux + sqrt( aux.^2 + R^2 - ( x(:,1).^2 +...
	 x(:,2).^2 )*ones(1,length(p_I)) ) ).^(-2*s) )*w_I;
ML = (areal*2*pi/s).*reshape( phi*weight , 3 , 3);
end
\end{Verbatim}
Recall the parametrization $\chi_{\ell}(\hat{x}) = B_{\ell}\hat{x} + x^{(1)}_{\ell}$, so that \verb+Bl+, \verb+x0+ and \verb+y0+ satisfy $\verb+Bl+ = B_{\ell}$ and $\left( \begin{array}{c} \verb+x0+ \\ \verb+y0+ \end{array} \right) = x^{(1)}_{\ell} .$ Moreover, \verb+s+ is $s$, \verb+areal+ is $|T_{ \ell}|$, \verb+p_I+ contains the quadrature points in the interval $[0,1]$, so that $2\pi\verb+p_I+(q) = \theta_q$, $\verb+w_I+(q) = W_q$, \verb+p_T+ contains 12 quadrature points over $\hat{T}$, stored in \verb+data.mat+ as \verb+p_T_12+ (see Appendix \ref*{sec:points}) and \verb+phi+ is the matrix $\Phi$, that is pre-computed and stored in \verb+data.mat+ as \verb+cphi+ (see Appendix \ref{sec:compdata}).

The output \verb+ML+ satisfies $\verb+ML+ \approx  2J_{\ell}$.    

%\bibliography{bib}{}
%\bibliographystyle{plain}

%%%%%%%%%%%%%%%%%%%%%%%%%%%%%%%%%%%%%%%%%%%%%%%%%%%%%%%%%%%%%%%%%%%%%%%%%%%%%%%%
\section{Two auxiliary functions}
\label{sec:fun}

%%%%%%%%%%%%%%%%%%%%%%%%%%%%%%%%%%%%%%%%%%%%%%%%%%%%%%%%%%%%%%%%%%%%%%%%%%%%%%%%%
The main code uses two functions that have not been outlined yet. Here we show them in detail.

The function \verb+setdiff_+ takes as input two vectors \verb+A+ and \verb+B+, such 
that \verb+A+ contains consecutive positive integers, ordered low to high, \verb+B+ contains 
positive integers and is such that $\verb+length(B)+ \leq \verb+length(A)+$ and $\verb+max(B)+ \leq \verb+max(A)+$. 
The function computes the set difference $\verb+A+ \setminus \verb+B+$, taking advantage 
of the pre-condition. 

\begin{Verbatim}[fontsize=\small]
function e = setdiff_( A , B )
e = A;
b = B - A(1) + 1;
b( b<1 )=[];
e(b) = [];
end   
\end{Verbatim}

%
%The function \verb+ball+ takes as inputs a positive real number \verb+R+ and a positive integer \verb+q+, and returns \verb+q+ equally spaced points over the circumference with radius \verb+R+ and centered in $(0,0)$. These points are saved as \verb+balln+. The function also returns the boundary segments that join consecutive points in the circumference, stored in \verb+balle+.     
%
%\begin{Verbatim}[fontsize=\small]
%function [balln balle] = ball( R , q )
%balln = R.*[ cos( 2*pi.*(0:1/q:1-1/q)); sin( 2*pi.*(0:1/q:1-1/q) )]';
%balle = [ 1:q ; [2:q 1] ]';
%end
%\end{Verbatim}

On the other hand, the function \verb+fquad+ calculates the entries of the right hand side vector in \eqref{eq:sistema}. 
Taking as input $\verb+areal+:=|T_{\ell}|$, the vectors \verb+xl+ and \verb+yl+, that contain the $x$ and $y$ coordinates of the vertices respectively, and a function $f$, \verb+fquad+ returns a vector in $\R^3$ array such that
$$ \verb+fquad+_k \approx \int_{T_{\ell}} f \, \varphi_{i_k}.$$
Here, for  $k \in \{1,2,3\}$, $i_k$ denotes the index of the $k$-th vertex of $T_{\ell}$ and $\phii_{i_k}$ the basis function corresponding to it.        
\begin{Verbatim}[fontsize=\small]
function VL = fquad( areal, xl , yl , f )
VL = zeros(3,1);
xmid = [(xl(2)+xl(3))/2, (xl(1)+xl(3))/2, (xl(1)+xl(2))/2];
ymid = [(yl(2)+yl(3))/2, (yl(1)+yl(3))/2, (yl(1)+yl(2))/2];
for i=1:3
    for j=1:3
        if j~=i
            VL(i) = VL(i) + areal/6 * f(xmid(j), ymid(j));
        end
    end
end
end
\end{Verbatim}

%%%%%%%%%%%%%%%%%%%%%%%%%%%%%%%%%%%%%%%%%%%%%%%%%%%%%%%%%%%%%%%%%%%%%%%%%%%%%%
\section{Auxiliary data}
\label{sec:data}
%%%%%%%%%%%%%%%%%%%%%%%%%%%%%%%%%%%%%%%%%%%%%%%%%%%%%%%%%%%%%%%%%%%%%%%%%%%%%%%
In order to perform the necessary calculations efficiently, along the execution the code makes use of pre-computed data, stored in \verb+data.mat+. 
Here we describe the variables provided by this file.
It is convenient to clarify that all the {\it MATLAB}\textsuperscript{\textregistered} code showed in this section does not belong to the program itself. It is included with an illustrative purpose.

\subsection{Quadrature points and weights: \texorpdfstring{\texttt{p\_cube}, \texttt{p\_T}, \texttt{p\_T\_comp}, \texttt{p\_I} and \texttt{w\_I}}{}}
%\subsection{Quadrature points and weights:} \texttt{p\_cube}, \texttt{xi}, \verb+xi_comp+, \texttt{xg} and \texttt{wg}\\
\label{sec:points}

We list the quadrature points used in all the quadrature rules and their respective weights. 

The matrix \verb+p_cube+ is used as input on functions \verb+vertex_quad+ and \verb+edge_quad+, and contains $27$ quadrature points over $[0,1]^3$. 
\begin{Verbatim}[fontsize=\small]
p_cube =

    0.1127    0.1127    0.1127
    0.1127    0.1127    0.5000
    0.1127    0.1127    0.8873
    0.1127    0.5000    0.1127
    0.1127    0.5000    0.5000
    0.1127    0.5000    0.8873
    0.1127    0.8873    0.1127
    0.1127    0.8873    0.5000
    0.1127    0.8873    0.8873
    0.5000    0.1127    0.1127
    0.5000    0.1127    0.5000
    0.5000    0.1127    0.8873
    0.5000    0.5000    0.1127
    0.5000    0.5000    0.5000
    0.5000    0.5000    0.8873
    0.5000    0.8873    0.1127
    0.5000    0.8873    0.5000
    0.5000    0.8873    0.8873
    0.8873    0.1127    0.1127
    0.8873    0.1127    0.5000
    0.8873    0.1127    0.8873
    0.8873    0.5000    0.1127
    0.8873    0.5000    0.5000
    0.8873    0.5000    0.8873
    0.8873    0.8873    0.1127
    0.8873    0.8873    0.5000
    0.8873    0.8873    0.8873
\end{Verbatim}

Over $\hat{T}$, we use two different quadrature rules, with $6$ and $12$ points. The set of nodes \verb+p_T_6+ is used to compute the non-touching element case and \verb+p_T_12+ as an input on \verb+comp_quad+.   
\begin{Verbatim}[fontsize=\small]
p_T_6 =

    0.5541    0.4459
    0.5541    0.1081
    0.8919    0.4459
    0.9084    0.0916
    0.9084    0.8168
    0.1832    0.0916

p_T_12 =

    0.7507    0.2493
    0.7507    0.5014
    0.4986    0.2493
    0.9369    0.0631
    0.9369    0.8738
    0.1262    0.0631
    0.6896    0.6365
    0.3635    0.0531
    0.9469    0.3104
    0.3635    0.3104
    0.6896    0.0531
    0.9469    0.6365
\end{Verbatim}

The $9\times 1$ array \verb+p_I+ contains the quadrature points over $[0,1]$, and \verb+w_I+ is a $9 \times 1$ array that contains their respective weights. These variables are used as input on \verb+comp_quad+. The set of nodes \verb+p_I+ is also employed in \verb+triangle_quad+. 
%   
%\begin{Verbatim}[fontsize=\small]
%p_I =
%
    %0.5000
    %0.0820
    %0.9180
    %0.0159
    %0.9841
    %0.3379
    %0.6621    
    %0.8067 
    %0.1933
%
%w_I =
%
    %0.1651
    %0.0903
    %0.0903
    %0.0406
    %0.0406
    %0.1562
    %0.1562
    %0.1303
    %0.1303
%\end{Verbatim}
\begin{Verbatim}[fontsize=\small]
p_I =                 w_I =
							
    0.5000                0.1651
    0.0820                0.0903
    0.9180                0.0903
    0.0159                0.0406
    0.9841                0.0406
    0.3379                0.1562
    0.6621                0.1562
    0.8067                0.1303 
    0.1933                0.1303
\end{Verbatim}

\subsection{Auxiliary variables to compute non-touching elements case: \texorpdfstring{\texttt{phiA}, \texttt{phiB} and \texttt{phiD}}{}} 
\label{sec:nontouchingdata}
The variables \verb+phiA+, \verb+phiB+ and \verb+phiD+ play the role of $\Phi^A$ , $\Phi^B$ and $\Phi^D$ (defined in Appendix \ref{sec:nontouching}), respectively. We expose below the code used to set up these variables. 
We use the lists \verb+p_T_6+ and \verb+w_T_6+ of quadrature points and weights in $\hat{T}$ defined in Appendix \ref{sec:points}: 
%We use \verb+xi+ defined as before (the quadrature points in $\hat{T}$) and \verb+w+ (defined below) as the weights associated to those quadrature points: 
\begin{Verbatim}[fontsize=\small]
w_T_6 = zeros(6,1);
w_T_6(1) = 0.1117;
w_T_6(2) = w_T_6(1);
w_T_6(3) = w_T_6(1);
w_T_6(4) = 0.0550;
w_T_6(5) = w_T_6(4);
w_T_6(6) = w_T_6(4);

local = cell(1,6);
local{1} = @(x,y) 1-x;
local{2} = @(x,y) x-y;
local{3} = @(x,y) y;
local{4} = @(x,y) -(1-x);
local{5} = @(x,y) -(x-y);
local{6} = @(x,y) -y;

mat_loc = zeros(6);
for i = 1:6
    for j = 1:6
        mat_loc(i,j) = local{i}(p_T_6(j,1),p_T_6(j,2));
    end
end

W = w_T_6*(w_T_6');
M_aux = zeros(18);
N_aux = zeros(18);
L_aux = zeros(18);

phiB = zeros(9,36);
phiA = zeros(9,36);
phiD = zeros(9,36);

for i=1:3
    for j=1:3
        for k = 1:6
            for q=1:6
                M_aux( q + 6*(i-1) , k + 6*(j-1) ) =...
                W(q,k)*mat_loc(i,q)*mat_loc(j+3,k);
                N_aux( q + 6*(i-1) , k + 6*(j-1) ) =...
                W(q,k)*mat_loc(i,q)*mat_loc(j,q);
                L_aux( q + 6*(i-1) , k + 6*(j-1) ) =...
                W(q,k)*mat_loc(i+3,k)*mat_loc(j+3,k);
            end
        end
    end
end

for i=1:9
    [im jm] = ind2sub([3 3] , i);
    im = 6*(im - 1) + 1;
    jm = 6*(jm - 1) + 1;
    phiB(i,:) = reshape( M_aux( im:im+5 , jm:jm+5 ) , 1 , [] );
    phiA(i,:) = reshape( N_aux( im:im+5 , jm:jm+5 ) , 1 , [] );
    phiD(i,:) = reshape( L_aux( im:im+5 , jm:jm+5 ) , 1 , [] );
end
\end{Verbatim}

\subsection{Auxiliary variables to compute vertex-touching elements case: \texorpdfstring{\texttt{vpsi1} and \texttt{vpsi2}}{}} 
\label{sec:vertexdata} 
The variables \verb+vpsi1+ and \verb+vpsi2+ are used as arguments of the 
function \verb+vertex_quad+ and play the role of the matrices 
$\Psi^1$ and $\Psi^2$ defined in Appendix \ref{sec:vertex}. 
Below we show the code used to initialize these variables. 

First we define a variable \verb+w_cube+ that lists the weights associated with each quadrature point stored in \verb+p_cube+:

\begin{Verbatim}[fontsize=\small]
w_cube =

    0.0214
    0.0343
    0.0214
    0.0343
    0.0549
    0.0343
    0.0214
    0.0343
    0.0214
    0.0343
    0.0549
    0.0343
    0.0549
    0.0878
    0.0549
    0.0343
    0.0549
    0.0343
    0.0214
    0.0343
    0.0214
    0.0343
    0.0549
    0.0343
    0.0214
    0.0343
    0.0214
\end{Verbatim}

The following lines generate \verb+vpsi1+ and \verb+vpsi2+: 
\begin{Verbatim}[fontsize=\small]
psi_D1 = cell(5,1);
psi_D1{1} = @(x,y,z) y-1;
psi_D1{2} = @(x,y,z) 1-x;
psi_D1{3} = @(x,y,z) x;
psi_D1{4} = @(x,y,z) -y.*(1-z);
psi_D1{5} = @(x,y,z) -y.*z;

psi_D2 = cell(5,1);
psi_D2{1} = @(x,y,z) -(y-1);
psi_D2{2} = @(x,y,z) y.*(1-z);
psi_D2{3} = @(x,y,z) y.*z;
psi_D2{4} = @(x,y,z) -(1-x);
psi_D2{5} = @(x,y,z) -x;

vpsi1 = zeros(25,27);
vpsi2 = zeros(25,27);

for i = 1:5
    for j = 1:5

        f1 = @(x,y,z) psi_D1{i}(x,y,z).*psi_D1{j}(x,y,z).*y;
        f2 = @(x,y,z) psi_D2{i}(x,y,z).*psi_D2{j}(x,y,z).*y;
        
        vpsi1( sub2ind([5 5], i , j) , : ) =...
        ( f1( p_cube(:,1) ,p_cube(:,2) , p_cube(:,3)) ).*w_cube;
        vpsi2( sub2ind([5 5], i , j) , : ) =...
        ( f2( p_cube(:,1) , p_cube(:,2) , p_cube(:,3)) ).*w_cube;
        
    end
end
\end{Verbatim}

\subsection{Auxiliary variables to compute edge-touching elements case: \texorpdfstring{\texttt{epsi1}, ..., \texttt{epsi5}}{}} 
\label{sec:edgedata} 
The variables \verb+epsi1+, ..., \verb+epsi5+ are used as input on the function \verb+edge_quad+ and play the role of $\Psi^1$, ..., $\Psi^5$ defined in Appendix \ref{sec:edge}, respectively. 
The code employed to set up these variables is exhibited below. 
We used the variable \verb+w_cube+ defined in the previous sub-section (containing weights associated to quadrature points stored in \verb+p_cube+):

\begin{Verbatim}[fontsize=\small]
psi_D1 = cell(3,1);
psi_D1{1} = @(x,y,z) -x.*y;
psi_D1{2} = @(x,y,z) x.*(1-z);
psi_D1{3} = @(x,y,z) x.*z;
psi_D1{4} = @(x,y,z) -x.*(1-y);

psi_D2 = cell(3,1);
psi_D2{1} = @(x,y,z) -x.*y.*z;
psi_D2{2} = @(x,y,z) -x.*(1-y);
psi_D2{3} = @(x,y,z) x;
psi_D2{4} = @(x,y,z) -x.*y.*(1-z);

psi_D3 = cell(3,1);
psi_D3{1} = @(x,y,z) x.*y;
psi_D3{2} = @(x,y,z) -x.*(1-y.*z);
psi_D3{3} = @(x,y,z) x.*(1-y);
psi_D3{4} = @(x,y,z) -x.*y.*z;

psi_D4 = cell(3,1);
psi_D4{1} = @(x,y,z) x.*y.*z;
psi_D4{2} = @(x,y,z) x.*(1-y);
psi_D4{3} = @(x,y,z) x.*y.*(1-z);
psi_D4{4} = @(x,y,z) -x;

psi_D5 = cell(3,1);
psi_D5{1} = @(x,y,z) x.*y.*z;
psi_D5{2} = @(x,y,z) -x.*(1-y);
psi_D5{3} = @(x,y,z) x.*(1-y.*z);
psi_D5{4} = @(x,y,z) -x.*y;

epsi1 = zeros(16,27);
epsi2 = zeros(16,27);
epsi3 = zeros(16,27);
epsi4 = zeros(16,27);
epsi5 = zeros(16,27);

for i = 1:4
    for j = 1:4
        
        f1 = @(x,y,z) psi_D1{i}(x,y,z).*psi_D1{j}(x,y,z) .*(x.^2);
        f2 = @(x,y,z) psi_D2{i}(x,y,z).*psi_D2{j}(x,y,z) .* (x.^2).*y;
        f3 = @(x,y,z) psi_D3{i}(x,y,z).*psi_D3{j}(x,y,z) .* (x.^2).*y;
        f4 = @(x,y,z) psi_D4{i}(x,y,z).*psi_D4{j}(x,y,z) .* (x.^2).*y;
        f5 = @(x,y,z) psi_D5{i}(x,y,z).*psi_D5{j}(x,y,z) .* (x.^2).*y;
        
        epsi1( sub2ind([4 4], i , j) , : ) =...
        ( f1( p_cube(:,1) , p_cube(:,2) , p_cube(:,3)) ).*w_cube;
        epsi2( sub2ind([4 4], i , j) , : ) =...
        ( f2( p_cube(:,1) , p_cube(:,2) , p_cube(:,3)) ).*w_cube;
        epsi3( sub2ind([4 4], i , j) , : ) =...
        ( f3( p_cube(:,1) , p_cube(:,2) , p_cube(:,3)) ).*w_cube;
        epsi4( sub2ind([4 4], i , j) , : ) =...
        ( f4( p_cube(:,1) , p_cube(:,2) , p_cube(:,3)) ).*w_cube;
        epsi5( sub2ind([4 4], i , j) , : ) =...
        ( f5( p_cube(:,1) , p_cube(:,2) , p_cube(:,3)) ).*w_cube;
        
    end
end
\end{Verbatim}

\subsection{Auxiliary variables to compute identical elements case: \texorpdfstring{\texttt{tpsi1}, \texttt{tpsi2} and \texttt{tpsi3}}{}} 
\label{sec:identicaldata}
Here, the variables \verb+tpsi1+, \verb+tpsi2+ and \verb+tpsi3+ are used as inputs on the function \verb+triangle_quad+ and play the role of the matrices $\Psi^1$, $\Psi^2 $ and $\Psi^3$, defined in Appendix \ref{sec:identical}, respectively. We describe the code used to set up these variables, where we use the quadrature data  \verb+p_I+ and \verb+w_I+ introduced in Appendix \ref{sec:points}:
\begin{Verbatim}[fontsize=\small]
lambda_D1 = cell(3,1);
lambda_D1{1} = @(z) -z;
lambda_D1{2} = @(z) -(1-z);
lambda_D1{3} = @(z) 1;

lambda_D2 = cell(3,1);
lambda_D2{1} = @(z) -1;
lambda_D2{2} = @(z) (1-z);
lambda_D2{3} = @(z) z;

lambda_D3 = cell(3,1);
lambda_D3{1} = @(z) z;
lambda_D3{2} = @(z) -1;
lambda_D3{3} = @(z) 1-z;

tpsi1 = zeros(9,9);
tpsi2 = zeros(9,9);
tpsi3 = zeros(9,9);

for i = 1:3
    for j = 1:3
        
        f1 = @(z) lambda_D1{i}(z).*lambda_D1{j}(z);
        f2 = @(z) lambda_D2{i}(z).*lambda_D2{j}(z);
        f3 = @(z) lambda_D3{i}(z).*lambda_D3{j}(z);
        
        tpsi1( sub2ind([3 3], i , j) , : ) =  f1( p_I ).*w_I;
        tpsi2( sub2ind([3 3], i , j) , : ) =  f2( p_I ).*w_I;
        tpsi3( sub2ind([3 3], i , j) , : ) =  f3( p_I ).*w_I;
        
    end
end
\end{Verbatim}

\subsection{Auxiliary variable to compute quadrature over complement: \texorpdfstring{\texttt{cphi}}{}}
\label{sec:compdata}
The matrix $\Phi$, defined in Appendix \ref{sec:comp}, is stored as the variable \verb+cphi+ and used as input on the function \verb+comp_quad+. Before explaining the code we employed to build it, we define the 12 by 1 array \verb+w_T_12+ as the set of weights associated to the quadrature points stored in \verb+p_T_12+:
\begin{Verbatim}[fontsize=\small]
w_T_12 =

    0.1168
    0.1168
    0.1168
    0.0508
    0.0508
    0.0508
    0.0829
    0.0829
    0.0829
    0.0829
    0.0829
    0.0829
\end{Verbatim}

Then, the following lines generate \verb+cphi+: 
\begin{Verbatim}[fontsize=\small]
local = cell(1,3);
local{1} = @(x,y) 1-x;
local{2} = @(x,y) x-y;
local{3} = @(x,y) y; 

cphi = zeros(9,12); 

for i = 1:3
    for j = 1:3
        
        f1 = @(z,y) local{i}(z,y).*local{j}(z,y);
        cphi( sub2ind([3 3], i , j) , : ) =...
        f1( p_T_12(:,1) , p_T_12(:,2) ).*w_T_12;
        
    end
end
\end{Verbatim}

%%%%%%%%%%%%%%%%%%%%%%%%%%%%%%%%%%%%%%%%%%%%%%%%%%%%%%%%%%%%%%%%%%%%%%%%%%%%%%
\section{Main Code}
\label{sec:maincode}
%%%%%%%%%%%%%%%%%%%%%%%%%%%%%%%%%%%%%%%%%%%%%%%%%%%%%%%%%%%%%%%%%%%%%%%%%%%%%%%
For the sake of the reader's convenience, we include here the main code described in Sections \ref{sec:codigo1} and \ref{sec:mainloop}. 
%
%\peligro{(2) Corregido en linea 59}

\begin{Verbatim}[fontsize=\small]
1 clc
2 s = 0.5;
3 f = @(x,y) 1;
4 cns = s*2^(-1+2*s)*gamma(1+s)/(pi*gamma(1-s));
5 load('data.mat');
6 nn = size(p,2); 
7 nt = size(t,1) 
8 uh = zeros(nn,1);
9 K  = zeros(nn,nn);
10 b  = zeros(nn,1);
11 % Compute areas
12 area = zeros(nt,1);
13 for i=1:nt
14     aux = p( : , t(i,:) );
15     area(i) = 0.5.*abs( det(  [ aux(:,1) - aux(:,3)...
				  aux(:,2) - aux(:,3)]  ) );
16 end
17 % Build patches data structure
18 deg = zeros(nn,1);
19 for i=1:nt
20     deg( t(i,:) ) = deg( t(i,:) ) + 1;
21 end
22 patches = cell(nn , 1);
23 for i=1:nn
24     patches{i} = zeros( 1 , deg(i) );
25 end
26 for i=1:nt
27     patches{ t(i,1) }(end - deg( t(i,1) ) + 1) = i;
28     patches{ t(i,2) }(end - deg( t(i,2) ) + 1) = i;
29     patches{ t(i,3) }(end - deg( t(i,3) ) + 1) = i;
30     deg( t(i,:) ) = deg( t(i,:) ) - 1;
31 end
32 % Preallocate auxiliary memory
33 vl = zeros(6,2);
34 vm = zeros(6*nt,2);
35 norms = zeros(36,nt);
36 ML = zeros(6,6,nt);
37 empty = zeros(nt,1);
38 aux_ind = reshape( repmat( 1:3:3*nt , 6 , 1 ) , [] , 1 );
39 empty_vtx = zeros(2,3*nt);
40 BBm = zeros(2,2*nt);
41 for l=1:nt-nt_aux % Main Loop
42     edge =  [ patches{t(l,1)} patches{t(l,2)} patches{t(l,3)} ];
43     [nonempty M N] = unique( edge , 'first' );
44     edge(M) = [];
45     vertex = setdiff(  nonempty , edge  );
46     ll = nt - l + 1 - sum( nonempty>=l );
47     edge( edge<=l ) = []; 
48     vertex( vertex<=l ) = []; 
49     empty( 1:ll ) = setdiff_( l:nt , nonempty ); 
50     empty_vtx(: , 1:3*ll) = p( : , t( empty(1:ll) , : )' );
51     nodl = t(l,:);
52     xl = p(1 , nodl); yl = p(2 , nodl);
53     Bl = [xl(2)-xl(1) yl(2)-yl(1); xl(3)-xl(2) yl(3)-yl(2)]'; 
54     b(nodl) = b(nodl) + fquad(area(l),xl,yl,f); 
55     K(nodl, nodl) = K(nodl, nodl)...
 + triangle_quad(Bl,s,tpsi1,tpsi2,tpsi3,area(l),p_I)...
 + comp_quad(Bl,xl(1),yl(1),s,cphi,R,area(l),p_I,w_I,p_T_12); 
56     BBm(:,1:2*ll) = reshape( [ empty_vtx( : , 2:3:3*ll )...
 -  empty_vtx( : , 1:3:3*ll ) ,  empty_vtx( : , 3:3:3*ll )...
 - empty_vtx( : , 2:3:3*ll ) ] , [] , 2)' ; 
57     vl = p_T_6*(Bl') + [ ones(6,1).*xl(1) ones(6,1).*yl(1) ]; 
58     vm(1:6*ll,:) = reshape(...
 permute(...
 reshape( p_T_6*BBm(:,1:2*ll) , [6 1 2 ll] ) , [1 4 3 2] ) , [ 6*ll 2 ] )...
 + empty_vtx(: , aux_ind(1:6*ll) )'; 
59     norms(:,1:ll) = reshape( pdist2(vl,vm(1:6*ll,:)), 36 , [] ).^(-2-2*s); 
60     ML(1:3,1:3,1:ll) =  reshape( phiA*norms(:,1:ll) , 3 , 3 , [] ); 
61     ML(1:3,4:6,1:ll) =  reshape( phiB*norms(:,1:ll) , 3 , 3 , [] );
62     ML(4:6,4:6,1:ll) =  reshape( phiD*norms(:,1:ll) , 3 , 3 , [] );
63     ML(4:6,1:3,1:ll) =  permute( ML(1:3,4:6,1:ll) , [2 1 3] );
64     % Assembling stiffness matrix
65     for m=1:ll
66         order = [nodl t( empty(m) , : )];
67         K(order,order) = K(order,order)...
		 + ( 8*area(empty(m))*area(l) ).*ML(1:6,1:6,m);
68     end
69     for m=vertex
70         nodm = t(m,:);
71         nod_com = intersect(nodl, nodm);
72         order = [nod_com nodl(nodl~=nod_com) nodm(nodm~=nod_com)];
73         K(order,order) = K(order,order)...
 + 2.*vertex_quad(nodl,nodm,nod_com,p,s,vpsi1,vpsi2,area(l),area(m),p_cube);
74     end
75     for m=edge
76         nodm = t(m,:);
77         nod_diff = [setdiff(nodl, nodm) setdiff(nodm, nodl)];
78         order = [ nodl( nodl~=nod_diff(1) ) nod_diff  ];
79         K(order,order) = K(order,order)...
 + 2.*edge_quad(...
nodl,nodm,nod_diff,p,s,epsi1,epsi2,epsi3,epsi4,epsi5,area(l),area(m),p_cube);
80     end
81 end
82 uh(nf) = ( K(nf,nf)\b(nf) )./cns; 
83 trimesh(t(1:nt - nt_aux , :), p(1,:),p(2,:),uh); 
\end{Verbatim}

%\lipsum[3]

\bibliography{bib_50lines}
\bibliographystyle{plain}
\end{document}